\documentclass[10pt,a4paper]{article}
\usepackage[utf8]{inputenc}
\usepackage{amsmath}
\usepackage{amsfonts}
\usepackage{amssymb}
\usepackage[colorinlistoftodos]{todonotes}
\usepackage[colorlinks=true]{hyperref}
\hypersetup{
    colorlinks=blue,%
    citecolor=red,%
    filecolor=black,%
    linkcolor=blue,%
    urlcolor=blue
}
\usepackage[margin=1.3in]{geometry}

\usepackage[amsmath,thmmarks,hyperref]{ntheorem}
{
	\theoremstyle{nonumberplain}
	\theoremheaderfont{\bfseries}
	\theorembodyfont{\normalfont}
	\theoremsymbol{\mbox{$\Box$}}
	\newtheorem{pf}{Proof.}
}

\numberwithin{equation}{section}
\def\R{\mathbb{R}}
\def\e{\epsilon}
\def\pa{\partial}
\newtheorem*{conjecture}{Conjecture.}
\newtheorem{thm}{Theorem}[section]
\newtheorem{lem}{Lemma}[section]
\newtheorem{rem}{Remark}[section]

\newtheorem{proposition}{\indent Proposition}[section]

\newtheorem{step}{Step}

\usepackage{upgreek}
\usepackage{graphicx}

 \allowdisplaybreaks
\title{article}
\begin{document}
	\title{\bf A Liouville type theorem of the linearly perturbed Paneitz equation on $S^3$ } 
	\date{}
	\author{\medskip Shihong Zhang    }

	\renewcommand{\thefootnote}{\fnsymbol{footnote}}
	\footnotetext[1]{S. Zhang: mg1821015@smail.nju.edu.cn }
	\maketitle
	
	{\noindent\small{\bf Abstract:}
		We prove a Liouville type theorem for the linearly perturbed Paneitz equation: For $\e>0$ small enough, if $u_\e$ is a positive smooth solution of
		$$P_{S^3} u_\e+\e u_\e=-u_\e^{-7} \qquad \mathrm{~~on~~}S^3,$$
		where $P_{S^3}$ is the Paneitz operator of the round metric $g_{S^3}$, then $u_\e$ is constant. This confirms a conjecture proposed by Fengbo Hang and Paul Yang in [ Int. Math. Res. Not. IMRN, 2020 (11) ].
		  }
	
	\medskip 

{{\bf $\mathbf{2020}$ MSC:} 35B06,35B20 (58J37,58J70)}

\medskip 

	
	\section{Introduction}
	The Yamabe problem has been completely resolved by a series of papers \cite{aubin,Schoen,Trudinger,Yamabe}. The associated conformally covariant operator of second order on a closed manifold $(M^n,g)$ is the conformal Laplacian
	\begin{equation}
	-\frac{4(n-1)}{n-2}\Delta_g+R_g, \label{equ1.1}
	\end{equation}
	where $n\geq 3$ and $R_g$ is the scalar curvature of $g$. In 1980's,
	S. Paneitz introduced fourth oder conformally covariant operator in \cite{Paneitz}. On three dimensional closed manifold\,\,$(M,g)$, such an operator, called Paneitz operator for convenience, is given by
	\[
	P_gu={\Delta}_g^2u-\mathrm{div}_g(\frac{5}{4}R_gg(\nabla_g u,\cdot)-4 \mathrm{Ric}_g(\nabla_g u,\cdot))-\frac{1}{2}Q_gu
	\]
	with its associated $Q$-curvature defined by
	
	\[
	Q_g=-\frac{1}{4}\Delta_g R_g-2|\mathrm{Ric}|_g^2+\frac{23}{32}R_g^2.
	\]
	The Paneitz operator is conformally covariant: If $u, \varphi \in C^\infty(M)$ and $u>0$, there holds
	$$P_{u^{-4}g}(\varphi)=u^7P_g(u \varphi).$$
	In particular, if $\varphi=1$, then the above equation becomes
	$$P_g(u)=-\frac{1}{2}Q_{u^{-4}g} u^{-7}.$$

	On $(S^3,g_{S^3})$, up to a positive constant multiple the above equation becomes
	
	\begin{equation}
	\Delta_{g_{S^3}}^2u+\frac{1}{2}\Delta_{g_{S^3}}u-\frac{15}{16}u=-u^{-7}.\label{equ1.3}
	\end{equation}
	The associated functional is defined by
	$$E(u)=\int_{S^3}u~ P_{S^3}u ~dV_{g_{S^3}}.$$
	
	By the inverse of the stereographic projection $I:\R^3\rightarrow S^3\verb|\|\{N\}$, where  $N$ is the north pole of $S^3$, we pull back equation (\ref{equ1.3})\,into $\R^3$:	
\begin{equation}\label{equ1.4}
{\Delta}^2v=-v^{-7}, \quad v>0 \qquad \mathrm{~~in~~} \R^3,
\end{equation}
where $v(x)=\left(\frac{1+|x|^2}{2}\right)^{1/2}(u\circ I)(x)$. 
Readers are referred to \cite{Choi&Xu,Xu} for the classification of smooth solutions to (\ref{equ1.4}). Some closely related topics are in \cite{NNPY,Hyder and Wei} etc.

In \cite{Xu and Yang}, X. Xu and P. Yang studied the constant $Q$-curvature problem on three manifolds under the positivity of Paneitz operator. However, $P_{S^3}$  has a negative eigenvalue. After that, M. Zhu and P. Yang in \cite{Yang&Zhu} proved the following sharp Sobolev inequality on $S^3$ via the rearrangement method:
$$
	\inf_{\varphi>0,\varphi\in H^2(S^3)}\|\varphi^{-1}\|^2_{L^6(S^3)}\int_{S^3}\varphi ~P_{S^3}\varphi ~dV_{g_{S^3}}\geq-\frac{15}{16}[\mathrm{Vol}(S^3)]^{4/3}.
$$
Recently, 
 F. Hang and P. Yang \cite{Hang&Yang} 
  proposed another rearrangement 
   method to prove the same inequality: 
If $P_{S^3}$ is replaced by $P_{S^3}+\e$ for some small constant $\e>0$, then we study the minimizing problem:
\begin{equation}\label{perturbed_problem}
\inf_{\varphi>0,\varphi\in H^2(S^3)} \|\varphi^{-1}\|_{L^6(S^3)}^2 \int_{S^3} \varphi~ (P_{S^3}\varphi+\e\varphi)~dV_{g_{S^3}}.
\end{equation}
It is not hard to check that $P_{S^3}+\e$ satisfies the $P^{+}$ condition; see  F. Hang and P. Yang \cite[Definition 1.1]{Hang&Yang1} for the definition of $P^+$ condition. Thus, the above perturbed minimizing problem  is well-defined. Next, they claimed that all minimizers of \eqref{perturbed_problem} must be constant. Furthermore, they  conjectured that all positive critical points of \eqref{perturbed_problem} must be constant, which is exactly  \cite[Conjecture 1.1]{Hang&Yang} as follows:

\begin{conjecture}[F. Hang and P. Yang \cite{Hang&Yang}]
If $\epsilon>0$ is a small constant and $u_\epsilon$ is a positive smooth solution of
\begin{equation}
P_{S^3}u_\epsilon+\epsilon u_\epsilon=-u_\epsilon^{-7} \qquad \mathrm{~~on~~} S^3,\label{equ1.5}
\end{equation}
then $u_\epsilon$ must be equal to a constant. 
\end{conjecture}
Finally, if the conjecture is true, then letting $\e \searrow 0$ we obtain the sharp Sobolev inequality.

The purpose of this paper is to give an affirmative answer to the above conjecture via the method of moving planes.

A closely relate topic is the linearly perturbed problem for  for the Yamabe equation, which have been studied by various authors. One particular case of 
these results in \cite{Druet,Druet&Hebey,Druet&Hebey&Robert} etc., is to study the follow perturbation problem: Let $(M,g)$ be an $n$-dimensional smooth closed manifold for $n \geq 3$ and $\e \in \R_+$,
\begin{equation*}
L_gu_{\epsilon}+\epsilon u_{\epsilon}=u_{\epsilon}^{\frac{n+2}{n-2}}  , \quad u>0 \qquad \mathrm{on} \quad (M,g).
\end{equation*}
O. Druet and E. Hebey in \cite{Druet&Hebey} gave some blow up examples on $S^n$ for $n\geq6$, under the hypothesis that

\[	\limsup\limits_{\epsilon \to 0}\| u_{\epsilon}\|_{H^1(g_{S^n})}\leq C.\]

In comparison with our case, indeed we first prove that the full set of positive smooth solutions of \eqref{equ1.5} is compact in the $C^\infty$ topology. Explicitly,

	\begin{thm}\label{thm1.1}
			If $u_{\epsilon}$ is a positive smooth solution of  (\ref{equ1.5}),
			then for  any $k \in \mathbb{N}$, there exist a uniform constant $\epsilon_0'>0$, and $A_0, C_0>0$ depending only $\epsilon_0'$ and $k$,  such that  $\forall~\epsilon\in(0,\epsilon_0')$,
		\begin{equation}
	\frac{1}{A_0}\leq u_{\epsilon}\leq A_0,   \qquad  \| u_{\epsilon}\|_{C^k(g_{S^3})}\leq C_0.
		\end{equation}
		\end{thm}
	
	   The moving plane method has been developed many years before, for instance, see \cite{Caffarelli&Gidas&Spruck,Chen&Li,Lin,Li&Zhu} etc. In this paper, we focus on dimension three, though some arguments still work well in higher dimensions.

		\begin{thm}\label{thm1.2}
		If $u_{\epsilon}$ is a positive smooth solution of  (\ref{equ1.5})  then there exists $0<\epsilon_0<\epsilon_0'$ with $\e_0'$ given in Theorem \ref{thm1.1},
		such that   $\forall\, \epsilon\in(0,\epsilon_0),\,u_{\epsilon}$ is radial symmetric with respect to  every critical point of $u_\e$.
	\end{thm}
	
	With the radial symmetry of $u_\e$ at hand, we are now in a position to prove our main result.
	\begin{thm}\label{thm1.3}
		If $u_{\epsilon}$ is a positive smooth solution of (\ref{equ1.5}),
		then  $u_{\epsilon}=(\frac{15}{16}-\epsilon)^{-\frac{1}{8}}, ~\forall \, \epsilon\in(0,\epsilon_0)$, where $\e_0$ is given in Theorem \ref{thm1.2}.			
	\end{thm}
	
	The organization of this paper is as follows. In Section \ref{sect2}, we establish the compactness of the full set of positive smooth solutions to \eqref{equ1.5}. In Section \ref{sec3}, through the stereographic projection, we derive the equivalence between PDE \eqref{equ1.5} and its associated integral equation in $\R^3$ (see \eqref{equb} below), which is crucial in the analysis of the possible singularity which comes from the Kelvin transform. In Section \ref{sec4}, we use the moving plane method to prove the radial symmetry of positive solutions of \eqref{equb} with respect to  every critical point. With these preparations, via the the Kazdan-Warner  condition, we complete the proof of the main Theorem \ref{thm1.3}  in Section \ref{sec5}.

	\section{Compactness: proof of Theorem \ref{thm1.1}}\label{sect2}
	
	This section is devoted to the proof of the compactness of the full set of positive solutions to \eqref{equ1.5}, i.e., Theorem \ref{thm1.1}.

We first recall some preliminary results in F. Hang and P. Yang  \cite{Hang&Yang, Hang&Yang1}.

 \begin{lem}[\protect{Hang-Yang \cite[Corollary 2.1]{Hang&Yang1}}]\label{lem2.2}
	If $u\in H^2(M,g)$ satisfies
	 \[\|u^{-1} \|_{L^6(M,g)}\leq C,\quad \| u\|_{H^2(M,g)}\leq C,\]  then $|u| \geq C(M,g)>0$.
	\end{lem}

	\begin{lem}[\protect{Hang-Yang \cite[Corollary 7.1]{Hang&Yang1}}]\label{lem dsa}
	If $u\in H^{2}(S^3) $ such that $u(p)=0$ for some$~p\in S^3$, then $E(u)\geq0$, and $E(u)=0$ if and only if $u=C_\ast G_{p}$ for some $C_\ast\in \R$, where $G_{p}$ is the Green's function of $P_{S^3}$ at $p$.

	\end{lem}
	
	The Kazdan-Warner type condition can be stated as follows.
	\begin{lem}[\protect{Hang-Yang \cite[Lemma 5.1]{Hang&Yang}}]\label{lem jo}
	If $0<Q\in C^\infty(S^3)$ and $u \in H^2(S^3)$ satisfies $P_{g_{S^3}}u=-Q\,u^{-7}$ on $S^3$, then
	\[\int_{S^3}\langle\nabla Q,\nabla x_{i}\rangle_{g_{S^3}} u^{-6}\,\,dV_{g_{S^3}}=0, \quad 1 \leq i \leq 4,\]
	where $x_i $\,are  the coordinate functions in $\R^4$.
\end{lem}
	
	\begin{lem}\label{lem j}
		If $u_{\epsilon}$ is a positive smooth solution of (\ref{equ1.5}) and there exists some positive constants $p$ and $C$  such that 
		\[
		\limsup\limits_{\epsilon \to 0}\| u_{\epsilon}\|_{L^p(S^3)}\leq C,
		\]
		then for any $ k\in \mathbb{N}$, there exist positive constants $A_0=A_0(C),C_0=C_0(C,k)$ such that
		\begin{equation}
		\frac{1}{A_0}<u_{\epsilon}<A_0, \qquad    \| u_{\epsilon}\|_{C^k(S^3)}\leq C_0.
		\end{equation}
	\end{lem}
\begin{pf}

	We multiply (\ref{equ1.5}) by $ u_{\epsilon}$ and integrate it over $S^3$ to show
	
	\begin{equation}
		\int_{S^3}(\Delta_{g_{S^3}}u_{\epsilon})^2-\frac{1}{2}\vert \nabla_{g_{S^3}}u_{\epsilon}\vert^2-
	(\frac{15}{16}-\epsilon)u_{\epsilon}^2 \,\,dV_{g_{S^3}}
	=-\int_{S^3}u_{\epsilon}^{-6}\,\,dV_{g_{S^3}}.\label{equ2.1}
	\end{equation}
	Thus
	\[
	\begin{aligned}
\int_{S^3}(\Delta_{g_{S^3}}u_{\epsilon})^2\,\,dV_{g_{S^3}}&\leq
\int_{S^3}\frac{1}{2}\vert \nabla_{g_{S^3}}u_{\epsilon}\vert^2+\frac{15}{16}u_{\epsilon}^2\,\,dV_{g_{S^3}}\\
&\leq\frac{1}{2}\int_{S^3}(\Delta_{g_{S^3}}u_{\epsilon})^2\,\,dV_{g_{S^3}}+C\int_{S^3}u_{\epsilon}^2\,\,dV_{g_{S^3}},	 \end{aligned}
\]
which implies that
$$\int_{S^3}(\Delta_{g_{S^3}}u_{\epsilon})^2\,\,dV_{g_{S^3}}\leq
C\int_{S^3}u_{\epsilon}^2\,\,dV_{g_{S^3}}.$$
By the standard elliptic theory and Sobolev embedding, we obtain
\begin{equation}
\| u_{\epsilon}\|_{C^{\frac{1}{2}}(S^3)}\leq C\| u_{\epsilon}\|_{H^2(S^3)}
\leq C\| u_{\epsilon}\|_{L^2(S^3)}.\label{equ2.2}
\end{equation}
If $p\geq2$, then  it follows from \eqref{equ2.2} and H\"older's inequality that
	\[
	\| u_{\epsilon}\|_{C^{\frac{1}{2}}(S^3)}\leq C.
	\]
If $0<p<2$, then
\[
\| u_{\epsilon}\|_{C^{\frac{1}{2}}(S^3)}\leq C(\int_{S^3}u_{\epsilon}^2\,\,dV_{g_{S^3}})^\frac{1}{2}\leq C(\int_{S^3}u_{\epsilon}^p\,\,dV_{g_{S^3}})^\frac{1}{2}\| u_{\epsilon}\|_{C^{\frac{1}{2}}(S^3)}^{\frac{2-p}{2}},
\]
which yields
$$\| u_{\epsilon}\|_{C^{\frac{1}{2}}(S^3)}\leq C \|u_\epsilon\|_{L^p(S^3)}\leq C. $$
In both cases, again by (\ref{equ2.2}) we obtain
\begin{equation}
\| u_{\epsilon}\|_{H^2(S^3)}\leq
C\| u_{\epsilon}\|_{L^2(S^3)}\leq
C\| u_{\epsilon}\|_{C^{\frac{1}{2}}(S^3)}\leq C.\label{equ2.3}
\end{equation}

We apply (\ref{equ2.1}) to estimate

\begin{equation*}
\int_{S^3}u_{\epsilon}^{-6}\,\,dV_{g_{S^3}}\leq
\int_{S^3}\frac{1}{2}\vert \nabla_{g_{S^3}}u_{\epsilon}\vert^2+
(\frac{15}{16}-\epsilon)u_{\epsilon}^2\,\,dV_{g_{S^3}}\leq C
\| u_{\epsilon}\|_{H^2(S^3)}^2\leq C.
\end{equation*}
Then Lemma \ref{lem2.2} shows that there exits a positive constant $C$ independent of $\epsilon$ such that
\begin{equation}\label{equ2.5}
u_{\epsilon}(p)\geq C, \quad \forall ~ \epsilon>0, \quad \forall ~ p\in S^3.
\end{equation}

Hence, it follows from (\ref{equ2.2}) and (\ref{equ2.5}) that $u_{\epsilon}$ is uniformly bounded from below and above. Notice that
	\[
	\Delta_{g_{S^3}}^2u_{\epsilon}+\frac{1}{2}\Delta_{g_{S^3}}u_{\epsilon}-(\frac{15}{16}-\epsilon)u_{\epsilon}=-u_{\epsilon}^{-7}.
	\]
By the standard elliptic estimates, we obtain
\[
 \| u_{\epsilon}\|_{C^4(S^3)}\leq C\quad \Longrightarrow \quad \| u_{\epsilon}\|_{C^k(S^3)}\leq C.
\]	
This completes the proof.	
\end{pf}

\begin{rem}
By the above analysis, if  $ \max\limits_{S^3}\, u_{\epsilon}\rightarrow +\infty$ , then $ \min\limits_{S^3}\, u_{\epsilon}\rightarrow 0$ at the same time, as $\epsilon \rightarrow 0$ .
\end{rem}

Indeed, we can prove that $\|u_\epsilon\|_{L^2(S^3)}$ is uniformly bounded.
\begin{lem}\label{lem bg}
	If $u_{\epsilon}$ is a positive smooth solution of (\ref{equ1.5}),  then there exists a uniform constant $C>0$ such that $\limsup\limits_{\epsilon \to 0}\| u_{\epsilon}\|_{L^2(S^3)}\leq C$.
\end{lem}
\begin{pf}
By contradiction, there exists a sequence of positive constants $\{\epsilon_i\}$ with $\epsilon_{i}\to 0$ such that $\,\| u_{\epsilon_{i}}\|_{L^2(S^3)}\to+\infty$ as $i \to \infty$. 
If we let 
$$\tilde{u}_{\epsilon_{i}}=\frac{u_{\epsilon_{i}}}{\| u_{\epsilon_{i}}\|_{L^2(S^3)}},$$
 then
 \begin{equation}\label{eq:blow-up_sol}
	P_{S^3}\tilde{u}_{\epsilon_{i}}+\epsilon_{i} \tilde{u}_{\epsilon_{i}}=-\| u_{\epsilon_{i}}\|_{L^2(S^3)}^{-8}\tilde{u}_{\epsilon_{i}}^{-7}.
\end{equation}
	
	A similar argument in Lemma\,\ref{lem j} shows
	\[\| \tilde{u}_{\epsilon_{i}}\|_{C^{\frac{1}{2}}(S^3)}\leq C\| \tilde{u}_{\epsilon_{i}}\|_{H^2(S^3)}
	\leq C\| \tilde{u}_{\epsilon_{i}}\|_{L^2(S^3)}\leq C.
	\]
	Then up to a subsequence, there hold $\tilde{u}_{\epsilon_{i}}\rightarrow\tilde{u}_{\infty} \geq 0$ in $C^{0}(S^3)$ with $\|\tilde u_\infty\|_{L^2(S^3)}=1$, $\tilde{u}_{\epsilon_{i}}\rightharpoonup\tilde{u}_{\infty}$ in $H^{2}(S^3)$ as $i \to \infty$.
	
	Our discussion is divided into two cases.
		
	\textbf{Case 1.}\,\,$\tilde{u}_{\infty}>0$.
	
	By \eqref{eq:blow-up_sol}, it is not hard to show that $\,\forall ~\phi\in C^{\infty}(S^3),$
	\[
	\int_{S^3}\Delta_{g_{S^3}}\tilde{u}_{\epsilon_{i}}\Delta_{g_{S^3}}\phi-\frac{1}{2}\langle\nabla_{g_{S^3}}\tilde{u}_{\epsilon_{i}},\nabla_{g_{S^3}}\phi\rangle -(\frac{15}{16}-\epsilon_{i})\tilde{u}_{\epsilon_{i}}\phi\,\,dV_{g_{S^3}}=-\| u_{\epsilon_{i}}\|_{L^2}^{-8}\int_{S^3}\tilde{u}_{\epsilon_{i}}^{-7}\phi\,\,dV_{g_{S^3}}.	\]
Letting $i \to \infty$, we have\,$\int_{S^3}\phi \,\,P_{S^3}\tilde{u}_{\infty}\,\,dV_{g_{S^3}}=0$. However, this together with $\phi=1$ implies
 \[\int_{S^3}\tilde{u}_{\infty}\,\,dV_{g_{S^3}}=0.\]
This is a contradiction.
	
	\textbf{Case 2.}\,\, There exists $p\,\in S^3$ such that $\tilde{u}_{\infty}(p)=0$.
	
	It follows from Lemma \ref{lem dsa} that  $E(\tilde{u}_{\infty})\geq 0$. On the other hand,
	\[
	\begin{aligned}
	&\int_{S^3}(\Delta_{g_{S^3}}\tilde{u}_{\epsilon_{i}})^2-\frac{1}{2}\vert \nabla_{g_{S^3}}\tilde{u}_{\epsilon_{i}}\vert^2-
	(\frac{15}{16}-\epsilon_{i})\tilde{u}_{\epsilon_{i}}^2\,\,dV_{g_{S^3}}= -\| u_{\epsilon_{i}}\|_{L^2}^{-8}\int_{S^3}\tilde{u}_{\epsilon_{i}}^{-6}\,\,dV_{g_{S^3}}<0.\\
\end{aligned}
		\]
 	By the lower semicontinuity of the norm $\|\cdot\|_{H^2(S^3)}$, we have
		\[
		E(\tilde{u}_{\infty})=\int_{S^3}\tilde{u}_{\infty} P_{S^3}\tilde{u}_{\infty}\,\,dV_{g_{S^3}}\leq \liminf_{i\to \infty}E(\tilde{u}_{\epsilon_{i}})\leq 0.
		\]
	This implies that $E(\tilde{u}_{\infty})=0$.  Furthermore, again by Lemma \ref{lem dsa} we conclude that $\tilde{u}_{\infty}=C_\ast G_{p}$ for some negative constant $C_\ast$.
		 However, the Kazdan-Warner type condition in Lemma \ref{lem jo} shows
		\[
		\begin{aligned}
			\int_{S^3}\langle\nabla_{g_{S^3}}(\tilde{u}_{\epsilon_{i}}^8\epsilon_{i}+\| u_{\epsilon_{i}}\|_{L^2}^{-8}),\nabla_{g_{S^3}}x_{4}\rangle \tilde{u}_{\epsilon_{i}}^{-6}\,\,dV_{g_{S^3}}=0
			\quad \Longrightarrow \quad \int_{S^3}x_{4} \tilde{u}_{\epsilon_{i}}^2 \,\,dV_{g_{S^3}}=0.
		\end{aligned}
	\]
	Letting $i \to \infty$ we obtain
		\[
	\int_{S^3}x_{4} \tilde{u}_\infty^2 \,\,dV_{g_{S^3}}=0.
	\]
	This is impossible, since $G_p(\cdot)=-\|\cdot\|_{S^3}/(8\pi)$ and 
	\begin{align*}
	 \int_{S^3}x_{4} \tilde{u}_\infty^2 \,\,dV_{g_{S^3}}=\frac{C_\ast^2}{(4\pi)^2} \int_{\R^3}\frac{|y|^2-1}{1+|y|^2}\left(\frac{2}{1+|y|^2}\right)^4 dy<0.
	 \end{align*}
\end{pf}

Therefore, we combine Lemmas \ref{lem j} and \ref{lem bg} to conclude Theorem \ref{thm1.1}.

\section{An associated integral equation}\label{sec3}
	
	For each $\e>0$, without loss of generality, we assume the south pole $S$ of $S^3$ is a critical point of $u_\e$ and  let
  $I_{\epsilon}:\R^3\rightarrow S^3\setminus \{N\}$ be the inverse of the stereographic projection, where $N$ is the north polar of $S^3$.

	\begin{proposition}
	If $u_\e$ is a positive smooth solution of \eqref{equ1.5} and let 
	\begin{equation}\label{def:v_e}
	v_{\epsilon}(x)=\left(\frac{1+\vert x\vert^2}{2}\right)^{\frac{1}{2}}(u_{\epsilon}\circ I_{\epsilon})(x),
	\end{equation}
	 then $v_\epsilon$ solves
	\begin{equation}
			{\Delta}^2v_{\epsilon}+\epsilon\left(\frac{2}{1+\vert x\vert^2}\right)^4v_{\epsilon}=-v_{\epsilon}^{-7}
			\quad \mathrm{~~in~~} \quad \R^3\label{equ3.1}
	\end{equation}
	if and only if 
	\begin{equation}
	v_{\epsilon}(x)=\int_{\R^3}\vert x-y\vert U_{\epsilon}(y)dy,\label{equb}
	\end{equation}
		where 
		$$U_{\epsilon}(x)=\frac{1}{8\pi}\left[\epsilon\left(\frac{2}{1+\vert x\vert^2}\right)^4 v_{\epsilon}+v_{\epsilon}^{-7}\right].$$
	\end{proposition}
	\begin{pf}
	
It follows from Lemma \ref{lem j} that  $U_{\epsilon}\in C^{\infty}(\R^3)$ and $U_{\epsilon}(y)=O (|y|^{-7})$ as $|y|\to \infty$.

For the sufficiency part, if $v_{\epsilon}(x)=\int_{\R^3}\vert x-y\vert U_{\epsilon}(y)dy$, then $\forall\,\, 0<|h|<1$,
\[
\begin{aligned}
  \bigg\vert \frac{v_{\e}(x+he_{i})-v_{\e}(x)}{h}\bigg\vert\leq& \int_{\R^3}\left| \frac{x+he_{i}-y\vert-\vert x-y\vert}{h}\right| U_{\epsilon}(y)dy\\
  =&\int_{\R^3}\frac{\vert h+2(x_{i}-y_{i})\vert}{|x+he_{i}-y|+| x-y| }U_{\epsilon}(y)dy\\
  \leq& \int_{\R^3}\frac{1+2\vert x-y\vert}{\vert x-y\vert}U_{\epsilon}(y)dy<+\infty.
\end{aligned}
\]
 So the Lebesgue dominated convergence theorem gives
\[
\frac{\partial v_{\e}}{\partial x_{i}}=\lim_{ h\to 0 }\frac{v_{\e}(x+he_{i})-v_{\e}(x)}{h}=\int_{\R^3}\frac{x_{i}-y_{i}}{\vert x-y\vert}U_{\epsilon}(y)dy.
\]
Through a similar argument, we further prove 
\[
\begin{aligned}
\frac{\partial^2 v_{\e}}{\partial x_{i}\partial x_{j}}=\int_{\R^3}\left[\frac{\delta_{ij}}{\vert x-y\vert}-\frac{(x_{i}-y_{j})(x_{j}-y_{j})}{\vert x-y\vert^3}\right]U_{\epsilon}(y)dy
\end{aligned}
\]
and then
$$\Delta v_{\e}(x)=2\int_{\R^3}\frac{U_{\epsilon}(y)}{\vert x-y\vert}dy.$$

Hence, the classic Newton potential theory shows
\[
{\Delta}^2v_{\epsilon}=-\epsilon\left(\frac{2}{1+\vert x\vert^2}\right)^4v_{\epsilon}-v_{\epsilon}^{-7}.
\]

However, the proof of the necessity part is very lengthy.  Somewhat inspired by Y. Choi- X. Xu \cite{Choi&Xu} and X. Xu \cite{Xu}, we first need to set up a series of preliminary estimates. 
	\end{pf}

\begin{lem}\label{lemh}
	If  $v_{\epsilon}$\ is a positive smooth solution of (\ref{equ3.1}) , then\, $\Delta v_{\epsilon}>0$ in $\R^3$.
\end{lem}
\begin{pf}
	If there exists $x_0\in\R^3$ such that $\Delta v_{\epsilon}(x_0)<0$, then we consider the spherical averages of 
	\[\bar{v}_{\epsilon}(r)=\frac{1}{\vert \partial B_{r}(x_0)\vert}\int_{\partial B_{r}(x_0)}v_{\epsilon} d \sigma,\quad
	\bar{w}_{\epsilon}(r)=\frac{1}{\vert \partial B_{r}(x_0)\vert}\int_{\partial B_{r}(x_0)}\Delta v_{\epsilon} d\sigma.
	\]
	
	We rewrite \eqref{equ3.1} as
	\[
	\Delta^{2} v_{\epsilon}=-\left(\frac{2}{1+\vert x\vert^2}\right)^{\frac{7}{2}}\big( (u_{\epsilon}\circ I_{\epsilon})^{-7}+\epsilon u_{\epsilon}\circ I_{\epsilon}\big)<0,
	\]
then we have
	\[
	\begin{aligned}
	\Delta \bar{v}_{\epsilon}=\bar{w}_{\epsilon}\quad \Longrightarrow \quad \bar{w}_{\epsilon}(0)=\Delta v_{\epsilon}(x_0)<0
	\end{aligned}
	\]
	and
	\[\Delta \bar{w}_{\epsilon}<0 \quad \mathrm{i.e.}\quad
	\frac{1}{r^2}\partial_{r}(r^2\partial_{r}\bar{w}_{\epsilon})<0.\]
	Also notice that $,\partial_{r}\bar{w}_{\epsilon}(0)=0$, then it is not hard to see that $\partial_{r}\bar{w}_{\epsilon}\leq0,\Delta \bar{v}_{\epsilon}(r)=\bar{w}_{\epsilon}(r)\leq\bar{w}_{\epsilon}(0)<0$. Moreover, we have 
	\begin{align*}
	\frac{1}{r^2}\partial_{r}(r^2\partial_{r}\bar{v}_{\epsilon})\leq \bar{w}_{\epsilon}(0), ~~\pa_r \bar v_\e(0)=0 \quad\Longrightarrow \quad 
	\partial_{r}\bar{v}_{\epsilon}(r)\leq\frac{\bar{w}_{\epsilon}(0)r}{3}.
	\end{align*}
	This again implies that
	$$
	\bar{v}_{\epsilon}(r)\leq \bar{v}_{\epsilon}(0)+\frac{\bar{w}_{\epsilon}(0)r^2}{6}=v_{\epsilon}(x_0)+\frac{\Delta v_{\epsilon}(x_0)r^2}{6}.
	$$
	Letting $r\to \infty$ we obtain a contradiction. 
	
	Finally, the strong maximum principle will imply $\Delta v_{\epsilon}>0$ in $\R^3$.
\end{pf}	 
\begin{lem}\label{lem:v_w}
	With the same assumption and the same notations as in Lemma \ref{lemh}, then
	 \[\partial_{r}\bar{v}_{\epsilon}(r)>0,\quad\partial_{r}^2\bar{v}_{\epsilon}(r)>0,\quad \partial_{r}^3\bar{v}_{\epsilon}(r)<0,\quad \bar{w}_{\epsilon}(r)>0, \quad \partial_{r}\bar{w}_{\epsilon}(r)<0\]
	 for $r>0$.
\end{lem}
	 \begin{pf}
	 	Again by \eqref{equ3.1} we have 
		$$\frac{1}{r^2}\partial_{r}(r^2\partial_{r}\bar{w}_{\epsilon}(r))<0.$$
Integrating the above inequality from $0$ to $r$ to show that  $\partial_{r}\bar{w}_{\epsilon}(r)<0$ for all $r>0$. By Lemma \ref{lemh}, we know
	 	 \[
	 	 \bar{w}_{\epsilon}(r)=\Delta \bar v_{\epsilon}(r)=\frac{1}{r^2}\partial_{r}(r^2\partial_{r}\bar{v}_{\epsilon}(r))>0.
	 	 \]
	 	This implies $\partial_{r}\bar{v}_{\epsilon}(r)>0$. 
		
		A direct calculation yields
		 $$\Delta^2\bar{v}_{\epsilon}=\frac{1}{r^4}\partial_{r}(r^4\partial_{r}^3\bar{v}_{\epsilon}(r))<0\quad \Longrightarrow \quad \partial_{r}^3\bar{v}_{\epsilon}(r)<0, ~~\forall ~ r>0.$$
	 	This shows $\partial_{r}^2\bar{v}_{\epsilon}(r)$ is decreasing in $(0,+\infty)$. We claim that $\pa_r^2 \bar v_\e(r)>0$ for all $r>0$. Otherwise, if there exists $r_{0}>0$ such that $\partial_{r}^2\bar{v}_{\epsilon}(r_{0})<0$, then
	 	\[ -\partial_{r}\bar{v}_{\epsilon}(r_{0})<\partial_{r}\bar{v}_{\epsilon}(r)-\partial_{r}\bar{v}_{\epsilon}(r_{0})\leq\partial_{r}^2\bar{v}_{\epsilon}(r_{0})(r-r_{0}).
	 	\]
	 	Letting  $r\to \infty$, we obtain a contradiction.
	 \end{pf}
	 
	 \begin{lem}\label{lemj}
	 There holds
	 	$\Delta v_{\epsilon}=2\int_{\R^3}\frac{1}{\vert x-y\vert} U_{\epsilon}(y)dy.$
	 \end{lem}
	 \begin{pf}
	 	For brevity, we set $h_{\epsilon}(x)=2\int_{\R^3}\frac{1}{\vert x-y\vert} U_{\epsilon}(y)dy$. We first 
	 claim that
	 \begin{equation}\label{est:h_e}
	  \lim_{ \vert x \vert\to \infty }\vert x\vert h_{\epsilon}(x)=2\alpha \quad\mathrm{with}\quad \alpha=\int_{\R^3}U_{\epsilon}(y)dy.
	  \end{equation}
	  
 For $\theta\in (0,1)$ and $|x|\gg1$, we split $\R^3=\cup_{i=1}^3 A_i$, where
 \begin{align*}
 A_1=&\{\vert x-y\vert<(1-\theta)\vert x\vert\}; \quad A_2=\{\vert y\vert<\theta \vert x\vert\}; \\
 A_3=&\{\vert x-y\vert\geq (1-\theta)\vert x\vert,\vert y\vert\geq \theta\vert x\vert\}.
 \end{align*}
 Then we have
	 	\begin{align*}
	 	\int_{\R^3}\frac{1}{\vert x-y\vert} U_{\epsilon}(y)dy
	 	=\sum_{i=1}^3 \int_{A_i}\frac{1}{\vert x-y\vert} U_{\epsilon}(y)dy.
		\end{align*}
		On $A_1$,  there holds $\theta |x|\leq |y|\leq (2-\theta)|x|$ and then
		$$ \int_{A_1}\frac{1}{\vert x-y\vert} U_{\epsilon}(y)dy\geq \frac{c_{\theta}}{\vert x\vert^7}\int_{A_1}\frac{1}{\vert x-y\vert}dy.$$
		On $A_2$, we estimate
		\begin{align*}
		&\int_{A_1}\frac{1}{\vert x-y\vert} U_{\epsilon}(y)dy\\
	 	\geq&\frac{1}{(1+\theta)\vert x\vert}(\int_{\R^3}-\int_{\vert y\vert>\theta\vert x\vert}U_{\epsilon}(y)dy)=\frac{\alpha}{(1+\theta)|x|}+\frac{c_\theta}{|x|^5}.
	 	\end{align*}
		Notice that $A_1\cup A_2 \subset B_{(1+\theta)|x|}(x)$, then $ B_{(1+\theta)|x|}(x)^c \subset A_3$. On $B_{(1+\theta)|x|}(x)^c$, there holds $|y|\geq |x-y|-|x|\geq (1-\frac{1}{1+\theta})|x-y|=\frac{\theta}{1+\theta}|x-y|$, then we estimate
		\begin{align*}
		\int_{A_3}\frac{1}{\vert x-y\vert} U_{\epsilon}(y)dy\geq& \int_{\{y; |x-y|<(1+\theta)|x|\}}\frac{1}{\vert x-y\vert} U_{\epsilon}(y)dy\\
	 	\geq&c_{\theta}\int_{\vert x-y\vert>(1+\theta)\vert x\vert}\frac{1}{\vert x-y\vert^8}
	 	=\frac{c_{\theta}}{\vert x\vert^5}.
		\end{align*}
		Combining the above estimates to show
		$$\int_{\R^3}\frac{1}{\vert x-y\vert} U_{\epsilon}(y)dy\geq \frac{\alpha}{(1+\theta)|x|}+\frac{c_\theta}{|x|^5}.$$
	 			A similar argument also yields
	 		$$	 		
				 		\int_{\R^3}\frac{1}{\vert x-y\vert} U_{\epsilon}(y)dy\leq \frac{{\alpha}}{(1-\theta)\vert x\vert}+\frac{c_{\theta}}{\vert x\vert^5}.
	 		$$
			Therefore, we conclude that
			$$\frac{2\alpha}{(1-\theta)}+\frac{c_\theta}{|x|^4}\leq |x|h_\e(x)\leq \frac{2\alpha}{(1+\theta)}+\frac{c_\theta}{|x|^4}.$$
			By letting $|x| \to \infty$ first and $\theta \searrow 0$ next, we prove the above claim.

		It follows from Lemma \ref{lemh} and \eqref{est:h_e} that $\Delta v_{\epsilon}-h_{\epsilon}>-C$, where $C$ is a positive constant independent of $\e$. Also notice that $\Delta (\Delta v_{\epsilon}-h_{\epsilon})=0$ in $\R^3$, then  there exists $C_\ast \in \R$ such that
	 	\[\Delta v_{\epsilon}=h_{\epsilon}+C_\ast=2\int_{\R^3}\frac{1}{\vert x-y\vert} U_{\epsilon}(y)dy+C_\ast.
	 	\]
	 	Lemma \ref{lem:v_w} states that $\partial_{r}\bar{w}_{\epsilon}(r)<0,\bar{w}_{\epsilon}(r)>0$ for all $r>0$, then 
		$$\lim\limits_{r\to\infty}\bar{w}_{\epsilon}(r)=a.$$
	 	If $a>0$, then 
		$$r^2\partial_{r}\bar{v}_{\epsilon}(r)=\int_{0}^{r}t^2\bar{w}_{\epsilon}(t)dt\geq\frac{a}{3}r^3\quad \Longrightarrow \quad \bar{v}_{\epsilon}(r)\geq \frac{a}{6}r^2, ~~\forall~ r>0.$$
		 However, it follows from  \eqref{def:v_e} and Theorem \ref{thm1.1} that $v_\e(x)=O(|x|)$ as $|x| \to \infty$. This yields a contradiction. 
		 
		 Hence, we conclude that
		  \[0=a=\lim\limits_{r\to\infty}\bar{w}_{\epsilon}(r)=\frac{1}{\vert \partial B_{r}(0)\vert}\int_{\partial B_{r}(0)}(h_{\epsilon}+C_\ast)d \sigma=C_\ast.\]
	 	This finishes the proof.
	 \end{pf}
	 \begin{lem} \label{lem:v_up_to_constant}
	 With a constant $\gamma_\e$, there holds
	 	$v_{\epsilon}(x)=\int_{\R^3}\vert x-y\vert U_{\epsilon}(y)dy+\gamma_{\epsilon}$.
	 \end{lem}
	 \begin{pf}
	 	For brevity, let $l_{\epsilon}(x)=\int_{\R^3}\vert x-y\vert U_{\epsilon}(y)dy$. Notice that
	 	 \[
	 	 \vert \nabla l_{\epsilon}(x)\vert\leq \int_{\R^3} U_{\epsilon}(y)dy= \alpha.
	 	 \]
	 	 This implies $l_{\epsilon}$ has at most linear growth at infinity, so does $v_\e-l_\e$. There also holds $\Delta(v_{\epsilon}-l_{\epsilon})=0$ in $\R^3$, this directly implies that $v_{\epsilon}(x)=l_{\epsilon}(x)+b_{i}x_{i}+\gamma_{\epsilon}$ for $b_i\in \R, 1 \leq i \leq 3$.
		 Moreover,
	 	 \[
	 	 \lim\limits_{\vert x\vert\to\infty}\frac{(v_{\epsilon}-l_{\epsilon})(x)}{|x|}=\frac{u_{\epsilon}(S)}{\sqrt{2}}-\frac{\alpha}{8\pi}.
	 	 \]
	 	 This indicates that $b_i=0, ~~1 \leq i \leq 3$.
	 \end{pf}
 \begin{lem}
 	$\vert \gamma_{\epsilon}\vert\leq C\epsilon. $
 \end{lem}
 \begin{pf}
 	 By Lemma \ref{lem:v_up_to_constant} we have
 	 \[\int_{B_R(0)}x\cdot\nabla v_{\epsilon} v_{\epsilon}^{-7} dx=\int_{B_{R}(0)}\int_{\R^3}\frac{\vert x\vert^2-x\cdot y
 	 }{\vert x-y \vert}U_{\epsilon}(y)v_{\epsilon}^{-7}(x)dydx,\]
	 
	 On one hand, an integration by parts gives
  \[
  LHS=-\frac{1}{6}\int_{B_R(0)}x\cdot\nabla v_{\epsilon} ^{-6}(x)dx=\frac{1}{2}\int_{B_{R}(0)} v_{\epsilon} ^{-6}dx-\frac{R}{6}\int_{\partial B_{R}(0)}v_{\epsilon}^{-6}d\sigma.
  \]
  Letting $R\to\infty$, we obtain  
  $$LHS=\frac{1}{2}\int_{\R^3} v_{\epsilon} ^{-6}dx.$$
  
  On the other hand,
  \[
  \begin{aligned}
  RHS=&\frac{1}{2}\int_{B_{R}(0)}\int_{\R^3}\frac{\vert x-y\vert^2-\vert y\vert^2+\vert x\vert^2
  }{\vert x-y \vert}U_{\epsilon}(y)\,v_{\epsilon}^{-7}(x)dy dx\\
=&\frac{1}{2}\int_{B_{R}(0)}(v_{\epsilon}-\gamma_{\epsilon})v_{\epsilon}^{-7}dx\\
&+\frac{1}{16\pi}\int_{B_{R}(0)}\int_{\R^3}
\frac{\vert x\vert^2-\vert y\vert^2
}{\vert x-y \vert}\left[\epsilon\left(\frac{2}{1+\vert x\vert^2}\right)^4 v_{\epsilon}(y)+v_{\epsilon}^{-7}(y)\right]v_{\epsilon}^{-7}(x)dydx.
  \end{aligned}
  \]
  Letting $R\to\infty$ and using the symmetry of the integral, we obtain 
  \[RHS=\frac{1}{2}\int_{\R^3}v_{\epsilon}^{-6}dx-\frac{\gamma_{\epsilon}}{2}\int_{\R^3}v_{\epsilon}^{-7}dx+\frac{\epsilon}{16\pi}\int_{\R^3}\int_{\R^3}
  \frac{\vert x\vert^2-\vert y\vert^2
  }{\vert x-y \vert}(\frac{2}{1+\vert x\vert^2})^4 v_{\epsilon}(y)v_{\epsilon}^{-7}(x)dxdy.\]

Therefore, we conclude that
\[
\begin{aligned}
\vert \gamma_{\epsilon}\vert&= \frac{\epsilon}{8\pi\int_{\R^3}v_{\epsilon}^{-7}dx}\left|\int_{\R^3}\int_{\R^3}
\frac{\vert x\vert^2-\vert y\vert^2
}{\vert x-y \vert}\left(\frac{2}{1+\vert x\vert^2}\right)^4 v_{\epsilon}(y)v_{\epsilon}^{-7}(x)dxdy\right|\\
&\leq C\epsilon \int_{\R^3}\int_{\R^3}
(\vert x\vert+\vert y\vert)
\left(\frac{2}{1+\vert x\vert^2}\right)^{7/2}\left (\frac{2}{1+\vert y\vert^2}\right)^{7/2}dxdy\\
&\leq C\epsilon.
\end{aligned}
\]
 \end{pf}

	 \begin{lem}[Pohoze\'av identity] \label{lem:Pohozeav}
	 
	 Let $Q \in C^\infty(\R^3)$ and $u$ be a positive smooth solution of 
	 	\[\Delta^2u+Q(x)u^{-7}=0 ,\,\, x\in \R^3,\]
	 	then
	 	\[
	 	-\frac{1}{6}\int_{B_{r}(0)}(x\cdot \nabla Q)u^{-6}dx=P(u,r)-\frac{1}{6}\int_{\partial B_{r}(0)}rQu^{-6}d\sigma,
	 	\]
	 	where
	 		\begin{align*}
	 	P(u,r)=\int_{\partial B_{r}(0)}\left[\frac{1}{2}r(\Delta u)^2-\frac{1}{2}u\frac{\partial\Delta u}{\partial r}	+\frac{1}{2}\Delta u\frac{\partial u}{\partial r}
	 	+r\frac{\partial u}{\partial r}\frac{\partial\Delta u}{\partial r}-\Delta u\frac{\partial}{\partial r}(r\frac{\partial u}{\partial r})\right]d\sigma.
	 	\end{align*}
	 \end{lem}
	 \begin{pf}
	 	It follows from a straightforward calculation.
	 \end{pf}
 
	 \begin{lem}\label{lemv}
	 There hold
	 	$$\lim\limits_{r\to\infty}P(v_{\epsilon},r)=4\pi\alpha\gamma_{\epsilon}$$
		and
		$$\,\lim\limits_{r\to\infty}r\int_{\partial B_{r}(0)}\left[\epsilon\left(\frac{2}{1+\vert x\vert^2}\right)^4v_{\epsilon}^8+1\right]v_{\epsilon}^{-6}d\sigma=0.$$
	 \end{lem}
 \begin{pf}
 	It suffices to prove the first assertion, since the second one is easy to show.
	
	It follows from Lemma \ref{lem:v_up_to_constant} that
 	\[
 	\begin{aligned}
 	v_{\epsilon}(x)&=\int_{\R^3}\vert x-y\vert U_{\epsilon}(y)dy+\gamma_{\epsilon}\\
 	&=\alpha\vert x\vert+\int_{\R^3}\big(\vert x-y\vert -\vert x\vert\big)U_{\epsilon}(y)dy+\gamma_{\epsilon}\\
 	&=\alpha\vert x\vert+\beta_{1}(x),
 	\end{aligned}
 	\]
and we estimate 
 \[
 \begin{aligned}
 | \beta_{1}(x)|&\leq\left| \int_{\R^3}\frac{\vert y\vert^2-2x\cdot y}{\vert x-y\vert +\vert x\vert}U_{\epsilon}(y)dy\right|+C\epsilon\\
 &\leq\int_{\R^3}\vert y\vert U_{\epsilon}(y)dy+C\epsilon\leq C.
 \end{aligned}
 \]

 Let $r=\vert x\vert$,
 \[
 \begin{aligned}
  \frac{\partial v_{\epsilon}}{\partial r}&= \int_{\R^3}\frac{x\cdot(x-y)}{\vert x-y\| x\vert}U_{\epsilon}(y)dy=\alpha+\int_{\R^3}\frac{x\cdot(x-y)-\vert x\| x-y \vert}{\vert x-y\| x\vert}U_{\epsilon}(y)dy\\
  &=\alpha+\frac{\beta_{2}(x)}{\vert x\vert},
  \end{aligned}
\]
we bound $\beta_{2}(x)$ by
\[
\begin{aligned}
|\beta_{2}(x)|&\leq\left| \int_{\R^3}\frac{\vert x-y\vert^2-\vert x\| x-y \vert+(x-y)\cdot y}{\vert x-y\vert }U_{\epsilon}(y)dy\right|\\
&\leq \int_{\R^3}\left| | x-y|-|x| \right| U_{\epsilon}(y)dy+\int_{\R^3}\vert y\vert U_{\epsilon}(y)dy\\
&\leq C.
\end{aligned}
\]

Notice that
\[
r\frac{\partial v_{\epsilon}}{\partial r}=\alpha\vert x\vert+\int_{\R^3}\big(\frac{x\cdot(x-y)}{\vert x-y\vert }-\vert x\vert \big)U_{\epsilon}(y)dy
\]
and then
\[
\begin{aligned}
\frac{\partial}{\partial r}\left(r\frac{\partial v_{\epsilon}}{\partial r}\right)&=\alpha+\int_{\R^3}\left(\frac{(2\vert x\vert^2-x\cdot y)\vert x-y\vert^2-(\vert x^2\vert-x\cdot y)^2}{\vert x-y\vert^3\vert x\vert}-1\right)U_{\epsilon}(y)dy\\
&=\alpha+\int_{\R^3}\left(\frac{(\vert x\vert^2-x\cdot y)(\vert y\vert^2-x\cdot y)}{\vert x-y\vert^3\vert x\vert}+\frac{\vert x\vert}{\vert x-y\vert}-1\right)U_{\epsilon}(y)dy\\
&=\alpha+\beta_{3}(x).
\end{aligned}
\]
We have
\[
\begin{aligned}
\left| |x|\beta_{3}(x)\right|&\leq\big\vert\int_{\R^3}\left[\frac{(\vert x\vert^2-x\cdot y)(\vert y\vert^2-x\cdot y)}{\vert x-y\vert^3}+\vert x\vert\frac{\vert x\vert-\vert x-y\vert}{\vert x-y\vert}\right]U_{\epsilon}(y)dy\big\vert\\
&\leq\int_{\R^3}\frac{\vert x\| y\vert}{\vert x-y\vert}U_{\epsilon}(y)dy+\int_{\R^3}\vert x\vert \frac{| y\cdot(x+x-y)|}{\vert x-y\vert(\vert x\vert+\vert x-y\vert)}U_{\epsilon}(y)dy\\
&\leq2\int_{\R^3}\frac{\vert x\| y\vert}{\vert x-y\vert}U_{\epsilon}(y)dy.
\end{aligned}
\]
Thus, we obtain
 \[\limsup_{\vert x\vert\to\infty}\big\| x\vert\beta_{3}(x)\big\vert\leq2\int_{\R^3}\vert y\vert U_{\epsilon}(y)dy\leq C.
 \]
 
 Similarly, by Lemma \ref{lemj} we have
 \[\begin{aligned}
 \Delta v_{\epsilon}&=2\int_{\R^3}\frac{1}{\vert x-y\vert} U_{\epsilon}(y)dy=\frac{2\alpha}{\vert x\vert}+2\int_{\R^3}\left(\frac{1}{\vert x-y\vert}-\frac{1}{\vert x\vert}\right)U_{\epsilon}(y)dy\\
 &=\frac{2\alpha}{\vert x\vert}+\beta_{4}(x).
 \end{aligned}
 \]
 We apply
 \[
 \begin{aligned}
 \big\| x\vert^2\beta_{4}(x)\big\vert&\leq2\int_{\R^3}\vert x\vert\frac{\big\| x\vert^2-\vert x-y\vert^2\big\vert}{| x-y\vert\big(\vert x\vert+\vert x-y\vert\big)}U_{\epsilon}(y)dy\\
 &\leq2\int_{\R^3}\frac{\vert x\vert \vert y\vert}{\vert x-y\vert}U_{\epsilon}(y)dy
 \end{aligned}
 \]
to show
 \[  \limsup_{\vert x\vert\to\infty}\left| |x|^2\beta_{4}(x)\right|\leq2\int_{\R^3}\vert y\vert U_{\epsilon}(y)dy\leq C.
 \]

Notice that
 \[
 \begin{aligned}
 \frac{\partial \Delta v_{\epsilon}}{\partial r}&=-2\int_{\R^3}\frac{x\cdot(x-y)}{\vert x-y\vert^3\vert x\vert}U_{\epsilon}(y)dy\\
 &=-\frac{2\alpha}{\vert x\vert^2}+2\int_{\R^3}\left(\frac{1}{\vert x\vert^2}-\frac{x\cdot(x-y)}{\vert x-y\vert^3\vert x\vert}\right)U_{\epsilon}(y)dy\\
 &=-\frac{2\alpha}{\vert x\vert^2}+\beta_{5}(x),
 \end{aligned}
  \]
Using
  \[
  \begin{aligned}
  \big\| x\vert^3\beta_{5}(x)\big\vert&\leq2\int_{\R^3}\left| |x|-\frac{\vert x\vert^2x\cdot(x-y)}{\vert x-y\vert^3}\right| U_{\epsilon}(y)dy\\
  &=2\int_{\R^3}\frac{\left| |x\| x-y\vert^3-\vert x\vert^2(\vert x\vert^2-x\cdot y)\right|}{\vert x-y\vert^3}U_{\epsilon}(y)dy\\
  &=2\int_{\R^3}\vert x\vert\frac{\left| |x-y\vert^3-\vert x\vert(\vert x-y\vert^2-\vert y\vert^2+x\cdot y)\right|}{\vert x-y\vert^3}U_{\epsilon}(y)dy\\
  &\leq 2\int_{\R^3}|x|\left(\frac{  |x-y\vert^2\big\| x-y\vert-\vert x\vert\big\vert}{\vert x-y\vert^3}+\frac{\vert x\vert \vert y\| x-y\vert}{\vert x-y\vert^3}\right)U_{\epsilon}(y)dy\\
  &\leq2\int_{\R^3}\left(\frac{\big(\vert x-y\vert+\vert x\vert\big)\vert x\| y\vert}{\vert x-y\vert\big(\vert x-y\vert+\vert x\vert\big)}+\frac{\vert x\vert^2\vert y\vert}{\vert x-y\vert^2}\right)U_{\epsilon}(y)dy,
  \end{aligned}
  \]
we obtain
  \[
   \limsup_{\vert x\vert\to\infty}\big\| x\vert^3\beta_{5}(x)\big\vert\leq4\int_{\R^3} \vert y\vert U_{\epsilon}(y)dy\leq C.
  \]
  
 We are now ready to calculate $P(v_{\epsilon},r)$:
  \[
  \begin{aligned}
  &\frac{1}{2}r(\Delta v_{\epsilon})^2-\frac{1}{2}v_{\epsilon}\frac{\partial\Delta v_{\epsilon}}{\partial r}	+\frac{1}{2}\Delta v_{\epsilon}\frac{\partial v_{\epsilon}}{\partial r}
  +r\frac{\partial v_{\epsilon}}{\partial r}\frac{\partial\Delta v_{\epsilon}}{\partial r}-\Delta v_{\epsilon}\frac{\partial}{\partial r}(r\frac{\partial v_{\epsilon}}{\partial r})\\
  =&\frac{1}{2}r(\frac{2\alpha}{r}+\beta_{4})^2-\frac{1}{2}(\alpha r+\beta_{1})(\frac{-2\alpha}{r^2}+\beta_{5})+\frac{1}{2}(\frac{2\alpha}{r}+\beta_{4})(\alpha+\frac{\beta_{2}}{r})\\
  &+(\alpha r+\beta_{2})(\frac{-2\alpha}{r^2}+\beta_{5})-(\frac{2\alpha}{r}+\beta_{4})(\alpha+\beta_{3})\\
  =&\alpha\big(\frac{3}{2}\beta_{4}-\frac{2\beta_{3}}{r}+\frac{\beta_{5}r}{2}+\frac{\beta_{1}-\beta_{2}}{r^2}\big)+(\frac{r}{2}\beta_{4}^2-\frac{1}{2}\beta_{1}\beta_{5}+\frac{\beta_{4}\beta_{2}}{r}+\beta_{2}\beta_{5}-\beta_{4}\beta_{3}).
  \end{aligned}
  \]
  Combining the above estimates of $\beta_{i}, 1 \leq i \leq 5$ to show
  \[\frac{r}{2}\beta_{4}^2-\frac{1}{2}\beta_{1}\beta_{5}+\frac{\beta_{4}\beta_{2}}{r}+\beta_{2}\beta_{5}-\beta_{4}\beta_{3}=O(\frac{1}{r^3}).
  \]
  We focus on the first term,
  \begin{align*}
  &\frac{3}{2}\beta_{4}-\frac{2\beta_{3}}{r}+\frac{\beta_{5}r}{2}+\frac{\beta_{1}-\beta_{2}}{r^2}\\
  =&3\int_{\R^3}\left(\frac{1}{\vert x-y\vert}-\frac{1}{\vert x\vert}\right)U_{\epsilon}(y)dy\\
  &-\frac{2}{\vert x\vert}\int_{\R^3}\left[\frac{(\vert x\vert^2-x\cdot y)(\vert y^2\vert-x\cdot y)}{\vert x-y\vert^3\vert x\vert}+\frac{\vert x\vert}{\vert x-y\vert}-1\right]U_{\epsilon}(y)dy\\
  &+\vert x\vert\int_{\R^3}\left(\frac{1}{\vert x\vert^2}-\frac{x\cdot(x-y)}{\vert x-y\vert^3\vert x\vert}\right)U_{\epsilon}(y)dy+\frac{\gamma_{\epsilon}}{r^2}\\
  &+\frac{1}{\vert x\vert^2}\int_{\R^3}\left[\big(\vert x-y\vert -\vert x\vert\big)-\frac{x\cdot(x-y)-\vert x\| x-y \vert}{\vert x-y\vert }\right]U_{\epsilon}(y)dy\\
  =&\frac{\gamma_{\epsilon}}{r^2}+\int_{\R^3}\left[\frac{\vert x\vert-\vert x-y\vert}{\vert x\| x-y\vert}-\frac{2(\vert x\vert^2-x\cdot y)(\vert y^2\vert-x\cdot y)}{\vert x-y\vert^3\vert x\vert^2}+\frac{1}{\vert x\vert}-\frac{x\cdot(x-y)}{\vert x-y\vert^3} \right.\\
  &\qquad\qquad~~\left.+\frac{\vert x-y\vert-\vert x\vert}{\vert x\vert^2} -\frac{x\cdot(x-y)-\vert x\| x-y \vert}{\vert x-y\| x\vert^2 }\right]U_{\epsilon}(y)dy\\
  =&\frac{\gamma_{\epsilon}}{r^2}+\int_{\R^3}\frac{(\vert y\vert^2-x\cdot y)}{\vert x-y\vert^3\vert x\vert^2}\left(\vert x-y\vert^2+\vert x\vert^2+2\vert y\vert^2-2x\cdot y-2\vert x-y\vert^2\right)U_{\epsilon}(y)dy\\
  =&\frac{\gamma_{\epsilon}}{r^2}+\int_{\R^3}\frac{\vert y\vert^2(\vert y\vert^2-x\cdot y)}{\vert x-y\vert^3\vert x\vert^2}U_{\epsilon}(y)dy.
   \end{align*}
    
    Also we estimate
   \[\big\vert\int_{\R^3}\frac{\vert y\vert^2(\vert y\vert^2-x\cdot y)}{\vert x-y\vert^3\vert x\vert^2}U_{\epsilon}(y)dy\big\vert\leq\frac{1}{|x|^2}\int_{\R^3}\frac{\vert y\vert^3}{\vert x-y\vert^2}U_{\epsilon}(y)dy=O(\frac{1}{r^4}).\]
    
    Therefore, putting these facts together we conclude that
    \[
    P(v_{\epsilon},r)=4\pi \alpha\gamma_{\epsilon}+O(\frac{1}{r}).
    \]
 \end{pf}
 \begin{lem}
 There holds
 	$\gamma_{\epsilon}=0$.
 \end{lem}
 \begin{pf}
 	Applying the Pohoza\'ev identity in Lemma \ref{lem:Pohozeav} with $Q=\epsilon (u_{\epsilon}\circ I_{\epsilon})^{8}+1$ and letting  $r \to \infty$, by Lemma \ref{lemv} we conclude that
 	\[
 	\lim_{r \to \infty}\int_{B_{r}(0)}\big(x\cdot\nabla Q\big)v_{\epsilon}^{-6} dx=-24\pi\alpha\gamma_{\epsilon}.
 	\]
		On the other hand, the Kazdan-Warner type condition gives
 \[
 \int_{S^3}\langle\nabla Q\circ I_\e^{-1},\nabla x_{4}\rangle_{g_{S^3}} u_{\epsilon}^{-6} dV_{g_{S^3}}=0.
 \]
Via stereographic projection, we obtain
 \[
 \int_{\R^3}\big(x\cdot\nabla Q\big)v_{\epsilon}^{-6} dx=0.
 \]
 
 Hence, $\gamma_\e=0$. 
 \end{pf}

 \section{Radial symmetry: Proof  of  the  Theorem 1.2}\label{sec4}
 
 We assume $S$ is a critical point of $u_{\e}$ and recall that
 $$v_{\epsilon}(x)=\left(\frac{1+\vert x\vert^2}{2}\right)^{\frac{1}{2}}(u_{\epsilon}\circ I_{\epsilon})(x)$$
as defined in \eqref{def:v_e}, then
 \begin{equation}
 \,\,\,\,\nabla v_{\epsilon}(0)=0.\label{equ3.2}
 \end{equation}
 We define the Kelvin transform of $v_{\epsilon}$ by 
 \begin{equation}
 v^{*}_{\epsilon}(x)=\vert x\vert v_{\epsilon}(\frac{x}{\vert x\vert ^2})=(\frac{1+\vert x\vert^2}{2})^{\frac{1}{2}}(u_{\epsilon}\circ I_{\epsilon})(\frac{x}{\vert x\vert^2}).\label{equ3.3}
 \end{equation}
Then $v_\e^\ast$ satisfies
 
 \begin{equation}
 {\Delta}^2v^{*}_{\epsilon}+\epsilon\left(\frac{2}{1+\vert x\vert^2}\right)^4v^{*}_{\epsilon}=-(v^{*}_{\epsilon})^{-7}
 , \qquad\mathrm{~~in~~} \quad\R^3\verb|\|\{{0}\}. \label{equ3.4}
 \end{equation} 
 Through (\ref*{equ3.2}) we have the following taylor expansion
 \begin{equation}\label{equ3.5}
 v^{*}_{\epsilon}(x)=a^{\epsilon}_0\vert x\vert+\frac{a^{\epsilon}_{ij}x_ix_j}{\vert x\vert^3}+O(\frac{1}{\vert x\vert^2}), \qquad\mathrm{~~for~~}\quad\vert x\vert\gg1,
 \end{equation}
 where $a^{\epsilon}_0=v_{\epsilon}(0),a^{\epsilon}_{ij}=\pa^2_{ij}v_{\epsilon}(0)/2$.
It follows from Theorem \ref{thm1.1} that there exist positive constants $A=A(A_0, C_0), C=C(A_0,C_0)$ independent of $\e$ such that
 \begin{equation}
 0<\frac{1}{A}<a^{\epsilon}_0<A,\,\,\,\vert a^{\epsilon}_{ij}\vert\leq C(C_0,A).\label{equ3.6}
 \end{equation}

	 	For $\lambda \leq 0$, we set $x^\lambda=(2\lambda-x_1,x_2,x_3), \varSigma_{\lambda}=\{ x \in \R^3\verb|||x_1<\lambda \}$ and  $T_{\lambda}=\{ x \in \R^3\verb|||x_1=\lambda \} $;
	 	 $v^{*}_{\epsilon,\lambda}(x)=v^{*}_{\epsilon}(x^{\lambda})$,
	 	 $w^{*}_{\epsilon,\lambda}(x)=v^{*}_{\epsilon,\lambda}(x)-v^{*}_{\epsilon}(x), x\in\varSigma_{\lambda}$.

		If $\lambda<0$, then $w^{*}_{\epsilon,\lambda}(x)$ has a possible singularity  at $d_{\lambda}=(2\lambda,0,0)$ and 
		\[
		\begin{aligned}
			\lim_{ x \to   d_{\lambda}}
		w^{*}_{\epsilon,\lambda}(x)&=\frac{u_{\epsilon}(N)}{\sqrt{2}}- \frac{\sqrt{1+4\lambda^2}}{\sqrt{2}}u_{\epsilon}\circ I_{\epsilon}(\frac{d_{\lambda}}{4\lambda^2}).\\
		\end{aligned}	\]
		Thus, we can assume $w^{*}_{\epsilon,\lambda} $ is a continuous function by defining \[w^{*}_{\epsilon,\lambda}(d_{\lambda})=\frac{u_{\epsilon}(N)}{\sqrt{2}}- \frac{\sqrt{1+4\lambda^2}}{\sqrt{2}}u_{\epsilon}\circ I_{\epsilon}(\frac{d_{\lambda}}{4\lambda^2}).\]

	Now we start with the method of  moving planes.

	\begin{lem}\label{lem3.1}
		There exist two constants $\lambda_0=\lambda_0(A_0,C_0)<0, R_0=R_0(A_0,C_0)>0$ independent of $\epsilon$  such that
		\[
		\Delta v^{*}_{\epsilon}(x)<\Delta v^{*}_{\epsilon}(x^{\lambda})  \qquad \forall~x \in \varSigma_{\lambda},\quad\lambda\leq\lambda_0,\quad\vert x^{\lambda}\vert>R_0.
		\]
	\end{lem}

	\begin{pf}
		For $|x|$ sufficiently large, a direct calculation yields
			\begin{equation}
					\begin{aligned}
				\Delta v^{*}_{\epsilon}(x)&=\frac{2}{\vert x\vert}v_{\epsilon}(\frac{x}{\vert x\vert^2})
				-\frac{4x\cdot(\nabla v_{\epsilon})(\frac{x}{\vert x\vert^2}) }{\vert x\vert^3}+
				\frac{1}{\vert x\vert^3}(\Delta v_{\epsilon})(\frac{x}{\vert x\vert^2})\\
				&=\frac{b_0^{\epsilon}}{\vert x\vert}+\frac{b_{ij}^{\epsilon}x_ix_j}{\vert x\vert^5}
				+O(\frac{1}{\vert x\vert^4}),
				\end{aligned}\label{a}
			\end{equation}
		where $b_0^{\epsilon}=2v_{\epsilon}(0)$, $b_{ij}^{\epsilon}$ only depends on $(v_{\epsilon})_{ij}(0)$. By Theorem \ref{thm1.1} we have
	\begin{equation}
		0<\frac{\sqrt{2}}{A_0}\leq b^{\epsilon}_0<\sqrt{2}A_0,\,\,\,\vert b^{\epsilon}_{ij}\vert\leq C(C_0,A).\label{b}
	\end{equation}
	The same argument can be applied to $\Delta v^{*}_{\epsilon}(x^{\lambda})$. The remaining part is standard, readers can refer to  \cite{Caffarelli&Gidas&Spruck,Lin} for details.
		\end{pf}
	
	\begin{lem}\label{lem3.2}
	There hold
	\begin{itemize}
		\item[(1)]$\Delta v^{*}_{\epsilon}\geq 0$ in  $\R^3\verb|\|\{0\};\,\,$
		\item[(2)]$
		\lim_{x\to0}\Delta v^{*}_{\epsilon}(x)=2\int_{\R^3}\vert y\vert^2 U_{\epsilon}(y)dy>0;
		$
		\item[(3)]$\forall~r>0$,\quad $\min\limits_{B_{r}\verb|\|\{0\}} 
		\Delta v^{*}_{\epsilon}=\min\limits_{\partial B_{r}} \Delta v^{*}_{\epsilon}$.
		\end{itemize}
	\end{lem}
	\begin{pf} Using (\ref{equb}), we calculate
	\[
	\begin{aligned}
		\nabla v_{\epsilon}(z)&=\int_{\R^3}\frac{z-y}{\vert z-y\vert}U_{\epsilon}(y)dy,\\
		\Delta v_{\epsilon}(z)&=2\int_{\R^3}\frac{1}{\vert z-y\vert}U_{\epsilon}(y)dy.
	\end{aligned}
	\]
	
	Letting $z=x/|x|^2$, by \eqref{a} we have
	\begin{align*}
			\Delta v^{*}_{\epsilon}(x)
			=&2\vert z\vert v_{\epsilon}(z)-4z\cdot (\nabla v_{\epsilon})(z)\vert z\vert+\vert z\vert^3\Delta v_{\epsilon}(z)\\
			=&2\int_{\R^3}|z|\left( \vert z-y\vert+\vert z\vert^2\frac{1}{\vert z-y\vert}-
		2\frac{z\cdot(z-y)}{\vert z-y \vert}\right)U_{\epsilon}(y)dy\\
		=&2\int_{\R^3}\frac{\vert z\vert}{\vert z-y\vert}\vert y\vert^2 U_{\epsilon}(y)dy\geq0.
			\end{align*}
		
	Thus,
	\[
	\lim_{x\to0}\Delta v^{*}_{\epsilon}(x)=\lim_{ \vert z\vert\to \infty}2\int_{\R^3}\frac{\vert z\vert}{\vert z-y\vert}\vert y\vert^2 U_{\epsilon}(y)dy=2\int_{\R^3}\vert y\vert^2 U_{\epsilon}(y)dy.
	\]
		
		It remains to show \textit{(3)}. For small $\delta >0$, we introduce
				$$\phi(x)=\frac{\delta}{\vert x\vert}+\Delta v^{*}_{\epsilon}-\min\limits_{\partial B_{r}} \Delta v^{*}_{\epsilon}.$$
				By \eqref{equ3.4}, $\phi$ satisfies
		\[
		\left\{\begin{aligned}
		\Delta \phi&<0\qquad\,\,\, \mathrm{in} \,\,\,B_{r}\verb|\|\{0\}, \\
		\phi&>0\qquad\,\,\, \mathrm{on}\,\,\, \partial B_{r}.
			\end{aligned}\right.	\]
	This together with \textit{(1)} implies that there exist a positive constant  $c_{\delta,r}$ such that $0<\rho<c_{\delta,r},\,\,\phi>0 \,on \,\,\partial B_{\rho}$. By the maximum principle, we have
		\[
		\phi>0\,\, in \,\,B_{r}\verb|\|B_{\rho},\,\,\,0<\rho<c_{\delta,r}.
		\]
By letting $\rho \to 0$ we obtain
		\[
		\phi\geq 0\,\, in \,\,B_{r}\verb|\| \{0\}.
		\]
		Next, letting $\delta \rightarrow 0$, we obtain the desired assertion.
	\end{pf}
	
	As before, we can assume  $\Delta v^{*}_{\epsilon}(x)$ is continuous function in $\R^3$ by defining $\Delta v^{*}_{\epsilon}(0)=2\int_{\R^3}\vert y\vert^2 U_{\epsilon}(y)dy$.

\vskip 8pt

	\begin{step}
		Let $w^{*}_{\epsilon,\lambda}(x)=v^{*}_{\epsilon}(x^{\lambda})-v^{*}_{\epsilon}(x)$, 
		then there exists a constant $\lambda_0'(A_0,C_0)<0$ such that
		
		\[
			\Delta w^{*}_{\epsilon,\lambda}>0 \quad\mathrm{and} \quad w^{*}_{\epsilon,\lambda}<0, \qquad \forall \,\,\lambda\leq\lambda_0',\quad x\in \varSigma_{\lambda}.
			\]
	\end{step}
	\vskip 8pt

		By Lemma \ref{lem3.1}, we know 
		$$\Delta w^{*}_{\epsilon,\lambda}>0,\quad \forall~x \in \varSigma_{\lambda},\quad\lambda\leq\lambda_0<0,\quad\vert x^{\lambda}\vert>R_0.$$ 
		  If $0<\vert x^{\lambda}\vert\leq R_0$, then by \eqref{a} and \eqref{b} we can choose $\lambda_0'<\lambda_0$ with $|\lambda_0'|$ such that
		  \[
		  \Delta v^{*}_{\epsilon}(x) < \frac{C}{|\lambda_0'|}\leq \min\limits_{\partial B_{R_0}} \Delta v^{*}_{\epsilon}\,\,\,,\forall~x\in \varSigma_{\lambda},\quad\lambda<\lambda_0'<\lambda_0.
		  \]
		  This together with Lemma \ref{lem3.2} \textit{(iii)} shows
		  \[
		  \Delta v^{*}_{\epsilon}(x) <  \min\limits_{\partial B_{R_0}} \Delta v^{*}_{\epsilon}\leq  \Delta v^{*}_{\epsilon}(x^{\lambda}) \quad \mathrm{i.e.} \quad\Delta w^{*}_{\epsilon,\lambda}>0,\quad  \forall~ x\in \varSigma_{\lambda}, \quad\lambda<\lambda_0'<\lambda_0.
		  \]
		  
		 Therefore, $\Delta w^{*}_{\epsilon,\lambda}>0,\forall~x\in \varSigma_{\lambda}$. Notice that $w^{*}_{\epsilon,\lambda}=0 \,\, on \,\,T_{\lambda}$, then it follows from \eqref{equ3.5} that
		  \[
		  \begin{aligned}
		  w^{*}_{\epsilon,\lambda}(x)&=a^{\epsilon}_0(\vert x^{\lambda}\vert-\vert x\vert)+O(\frac{1}{\vert x\vert})+O(\frac{1}{\vert x^{\lambda}\vert})\\
		  &=a^{\epsilon}_0\frac{4\lambda(\lambda-x_1)}{\vert x\vert+\vert x^{\lambda}\vert}+O(\frac{1}{\vert x^{\lambda}\vert}).
		  \end{aligned}
		  \]
		 Hence, for each fixed $\lambda<0$, there holds
		  \begin{equation}\label{est:w_at_infinity}
		  \liminf_{x\in\varSigma_{\lambda}, \vert x\vert\to \infty}w^{*}_{\epsilon,\lambda}(x)\leq0.
		  \end{equation}
	In particular,
		  \[
		  w^{*}_{\epsilon,\lambda}(d_\lambda)\leq \frac{A_0}{\sqrt{2}}- \frac{\sqrt{1+4\lambda^2}}{\sqrt{2}A_0}<0.\]
		  By the strong maximum principle, we obtain $w^{*}_{\epsilon,\lambda}(x)<0,\forall ~x\in \varSigma_{\lambda}$.

\begin{lem} \label{lem3.3}
	If  $w^{*}_{\epsilon,\lambda}(x)<0, \forall~ x\in \varSigma_{\lambda},\lambda \in [\lambda_0',0)$, then there exists $0<\epsilon_0=\e_0(A_0,C_0)<\epsilon_0'$ such that $$\Delta^2 w^{*}_{\epsilon,\lambda}<0,\quad\forall~ x\in \varSigma_{\lambda}\verb|\|\{d_{\lambda}\},\quad\forall~\epsilon  \in (0,\epsilon_0),$$where $A_0,C_0$ and $\e_0$ are given in Theorem \ref{thm1.1}.
\end{lem}	

	\begin{pf}
		The proof is complete by the following two claims.
	\vskip 8pt
	\textbf{Claim 1.} $\exists\,\, R_0''(A_0,C_0)>0,\,\,\,\Delta^2 w^{*}_{\epsilon,\lambda}<0,\forall~ x\in \varSigma_{\lambda},~\vert x\vert>R_0''>2\vert \lambda_{0}'\vert,\epsilon \in (0,\epsilon_0^{1})$\\
	\[
	\begin{aligned}
	\Delta^2 w^{*}_{\epsilon,\lambda}&=-(v^{*}_{\epsilon,\lambda})^{-7}+(v^{*}_{\epsilon})^{-7}+
	\epsilon\left(\frac{2}{1+\vert x\vert^2}\right)^4 v^{*}_{\epsilon}-	\epsilon\left(\frac{2}{1+\vert x^{\lambda}\vert^2}\right)^4 v^{*}_{\epsilon,\lambda}.\\
	&\leq
	\frac{(v^{*}_{\epsilon,\lambda})^{7}-(v^{*}_{\epsilon})^{7}}{(v^{*}_{\epsilon})^{7}(v^{*}_{\epsilon,\lambda})^{7}}-	\epsilon\left(\frac{2}{1+\vert x^{\lambda}\vert^2}\right)^4(v^{*}_{\epsilon,\lambda}-v^{*}_{\epsilon})\\
	&=\left(\frac{\sum_{i=0}^{6} (v^{*}_{\epsilon,\lambda})^i(v^{*}_{\epsilon})^{6-i}}{(v^{*}_{\epsilon})^{7}(v^{*}_{\epsilon,\lambda})^{7}}-	\epsilon(\frac{2}{1+\vert x^{\lambda}\vert^2})^4 \right)w^{*}_{\epsilon,\lambda}.
	\end{aligned}
	\]
By $\lambda \in [\lambda_0',0)$ and (\ref{equ3.5}), there exists a positive constant $R_0''=R_0''(\lambda_0')$ such that 
	\[
	\vert x\vert\sim \vert x^{\lambda}\vert,\,\,\, v^{*}_{\epsilon,\lambda}\sim \vert x\vert,v^{*}_{\epsilon}\sim \vert x\vert,\qquad \forall~\vert x\vert>R_0''.
	\]
	If we choose $0<\epsilon<\epsilon_0^{1}=\frac{C}{2C^{'}}$, then
	\[
	\Delta^2 w^{*}_{\epsilon,\lambda}\leq\left(\frac{C}{\vert x\vert^8}-\frac{C'\epsilon}{(1+\vert x \vert^2)^4}\right)w^{*}_{\epsilon,\lambda}<0,\,\,\,\,\,\forall \,\,\vert x\vert>R_0''.
	\]
	
	\vskip 8pt
	\textbf{Claim 2.}\,\,$\Delta^2 w^{*}_{\epsilon,\lambda}<0,\quad\forall~ x\in \varSigma_{\lambda},~\vert x\vert\leq R_0'',~\epsilon \in (0,\epsilon_0^{2})$.\\
	Notice that $\vert x^{\lambda}\vert<\vert x\vert\leq R_0''$ for $x\in \varSigma_{\lambda}$, by (\ref{equ3.3}) and Theorem \ref{thm1.1} we have 
	 \[
	C_2(A_0,R_0'',C_0)\geq v^{*}_{\epsilon}\geq C_1(A_0,R_0'',C_0).
	\] 
	If we choose $0<\epsilon<\frac{3(C_1)^{6}}{8(C_2)^{14}}=\epsilon_0^{2}$, then
	\[
	\Delta^2 w^{*}_{\epsilon,\lambda}\leq\left(\frac{6 (C_1)^{6}}{(C_2)^{14}}-16\epsilon\right)w^{*}_{\epsilon,\lambda}<0.
	\]
	
	Finally, letting $\epsilon_0=\min\{\epsilon_0^{2},\epsilon_0^{1}\}$, we obtain the desired assertion.
\end{pf}

		In the following, fix $\epsilon\in(0,\epsilon_0)$ and the constants $\delta_{\epsilon}^{1},\delta_{\epsilon}^{2}\cdots; R_{\epsilon}^{1}, R_{\epsilon}^{2} \cdots$ depend only on  $\epsilon$.

		We define
		 $$\bar{\lambda}_{\epsilon}=\sup\{\lambda\leq0\verb|| |\Delta w^{*}_{\epsilon,\mu}>0, w^{*}_{\epsilon,\mu}<0, \forall~x\in \varSigma_{\mu} , \mu\leq\lambda\},$$
		  then $0\geq\bar{\lambda}_{\epsilon}\geq \lambda_0'$.
		  
	\vskip 8pt
	\begin{step}
		If $\bar{\lambda}_{\epsilon}<0$, then $  \,\,w^{*}_{\epsilon,\bar{\lambda}_{\epsilon}}\equiv 0\,\, in \,\,\varSigma_{\bar{\lambda}_{\epsilon}}$.
	\end{step}
	\vskip 8pt
Now	we argue by contradiction, we assume $ w^{*}_{\epsilon,\bar{\lambda}_{\epsilon}}\not\equiv 0$. 
\begin{lem}\label{lem3.4}
		If $\lambda_0'<\bar{\lambda}_{\epsilon}<0$ and  $w^{*}_{\epsilon,\bar{\lambda}_{\epsilon}}\not\equiv 0$, then there exist $0<\delta_{\epsilon}^{1}\ll 1$  and $R^{1}_{\epsilon}\gg3\vert \lambda_0'\vert$ such that  such that
		\[
		\frac{\partial\Delta v^{*}_{\epsilon} }{\partial x_1}>0,\quad\forall~\vert x\vert>R^{1}_{\epsilon},~ \bar{\lambda}_{\epsilon}- \delta_{\epsilon}^{1}<x_1<\bar{\lambda}_{\epsilon}+ \delta_{\epsilon}^{1}<0.
		\]
\end{lem}
	\begin{pf}
		If we choose $R^{1}_{\epsilon}\gg1$, then $\vert x\vert\sim \vert x^{\bar{\lambda}_{\epsilon}}\vert$. By definition of $\bar{\lambda}_{\epsilon}$ and \eqref{est:w_at_infinity}, we have
		 \[
		 \begin{aligned}
& \Delta w^{*}_{\epsilon,\bar{\lambda}_{\epsilon}}\geq 0,\quad w^{*}_{\epsilon,\bar{\lambda}_{\epsilon}}\leq0,\,\,x\in\varSigma_{\bar{\lambda}_{\epsilon}},\\
& \Delta w^{*}_{\epsilon,\bar{\lambda}_{\epsilon}}=0,\quad w^{*}_{\epsilon,\bar{\lambda}_{\epsilon}}=0\,\, \,\, on \,\,T_{\bar{\lambda}_{\epsilon}},\\
 &\liminf_{x\in\varSigma_{\bar{\lambda}_{\epsilon}}, \vert x\vert\to \infty}w^{*}_{\epsilon,\bar{\lambda}_{\epsilon}}(x)\leq0.
		 \end{aligned} \]
		By the strong maximum principle, we have $w^{*}_{\epsilon,\bar{\lambda}_{\epsilon}}<0,\,\, x\in\varSigma_{\bar{\lambda}_{\epsilon}}$. Lemma \ref{lem3.3} implies that 
		
		\[\Delta^2 w^{*}_{\epsilon,\bar{\lambda}_{\epsilon}}<0,\quad \forall~ x\in \varSigma_{\bar{\lambda}_{\epsilon}}\verb|\|\{d_{\bar{\lambda}_{\epsilon}}\}.\]
		Again by the strong maximum principle, we obtain
		$$\Delta w^{*}_{\epsilon,\bar{\lambda}_{\epsilon}}> 0,\quad \forall~x\in\varSigma_{\bar{\lambda}_{\epsilon}}\verb|\|\{d_{\bar{\lambda}_{\epsilon}}\}.$$
		We apply a similar argument as in the proof of Lemma \ref{lem3.2} \textit{(3)}  to show		  \begin{equation}\label{equc}
		  	\min\limits_{B_{r}(d_{\bar{\lambda}_{\epsilon}})} 
		  	\Delta w^{*}_{\epsilon,\bar{\lambda}_{\epsilon}}=\min\limits_{\partial B_{r}(d_{\bar{\lambda}_{\epsilon}})} \Delta w^{*}_{\epsilon,\bar{\lambda}_{\epsilon}}>0,\,\,\,\,\forall\,\, 0<r<\vert \bar{\lambda}_{\epsilon}\vert,
		  \end{equation} 
		  thus $$\Delta w^{*}_{\epsilon,\bar{\lambda}_{\epsilon}}> 0,\quad \forall~x\in\varSigma_{\bar{\lambda}_{\epsilon}}.$$
		  Let $b_{\epsilon}=(\bar{\lambda}_{\epsilon},0,0), \bar{x}= x-b_{\epsilon}$, then
		  $\vert x\vert\sim \vert x^{\bar{\lambda}_{\epsilon}}\vert\sim\vert \bar{x}\vert$, if$\,\, \vert x\vert\gg \max\{\vert \lambda_{0}'\vert,R_{0}\}.$\\
		  
		  \textbf{Claim.} There exists a positive constant $\alpha_{\epsilon}$ such that 
		  \begin{equation}\label{est:w_ast}
		  \Delta w^{*}_{\epsilon,\bar{\lambda}_{\epsilon}}>\frac{\alpha_{\epsilon}(\bar{\lambda}_{\epsilon}-x_1)}{\vert \bar x\vert^3},\qquad  \forall~x\in\varSigma_{\bar{\lambda}_{\epsilon}}\setminus B_{R^{1}_{\epsilon}},
		  \end{equation}
		  where $B_{R^{1}_{\epsilon}}=B_{R^{1}_{\epsilon}}(0)$.
		 
		 The Hopf lemma implies
		  \[
		  \frac{\partial\Delta w^{*}_{\epsilon,\bar{\lambda}_{\epsilon}} }{\partial x_1}\leq c_{\epsilon}<0, ~~\forall~ x \in (B_{R^{1}_{\epsilon}+1}\verb|\|B_{R^{1}_{\epsilon}-1})\cap T_{\bar{\lambda}_{\epsilon}}.
		  \]

		  By the uniform continuity, we can choose  $0<\delta_{\epsilon}\ll 1 $ such that
		  \[
		  \frac{\partial\Delta w^{*}_{\epsilon,\bar{\lambda}_{\epsilon}} }{\partial x_1}\leq \frac{c_{\epsilon}}{2}<0,\quad \forall~ x \in D_1,
		  \]
		  where
		 \[
		 \begin{aligned}
		 D_1=&(B_{R^{1}_{\epsilon}+1}\verb|\|B_{R^{1}_{\epsilon}-1})\cap\{x\verb|||\bar{\lambda}_{\epsilon}- \delta_{\epsilon}<x_1<\bar{\lambda}_{\epsilon}+ \delta_{\epsilon}\},\\
		 D_2=&(\partial B_{R^{1}_{\epsilon}})\cap\{x\verb|||x_1\leq\bar{\lambda}_{\epsilon}- \delta_{\epsilon}\}.
		 \end{aligned}
		 \]

		  By direct calculation,
		  \[
		  \left| \frac{\partial}{\partial x_1}\frac{ (\bar{\lambda}_{\epsilon}-x_1)}{\vert \bar{x}\vert^3}\right|<\frac{C}{\vert\bar{x}\vert^3}.
		  \]
		  
		    We further choose $\alpha_{\epsilon}>0$ small enough such that 	
		    $$	  \frac{\partial}{\partial x_1}\left(\Delta w^{*}_{\epsilon,\bar{\lambda}_{\epsilon}}-\frac{\alpha_{\epsilon} (\bar{\lambda}_{\epsilon}-x_1)}{\vert \bar{x}\vert^3}\right)<0 \qquad \mathrm{for~~all~~} x \in D_1$$
		     and
		  \[
		  \Delta w^{*}_{\epsilon,\bar{\lambda}_{\epsilon}}-\frac{\alpha_{\epsilon} (\bar{\lambda}_{\epsilon}-x_1)}{\vert \bar{x}\vert^3}>0 \qquad \mathrm{for~~all~~} x \in D_2,
		  \]
		    then \[\Delta w^{*}_{\epsilon,\bar{\lambda}_{\epsilon}}-\frac{\alpha_{\epsilon} (\bar{\lambda}_{\epsilon}-x_1)}{\vert \bar{x}\vert^3}>0 \qquad \forall~x\in D_1\cap\partial B_{R^{1}_{\epsilon}}.\]
		    
		    Observe that
		      \[
		  \Delta\left(\Delta w^{*}_{\epsilon,\bar{\lambda}_{\epsilon}}-\frac{\alpha_{\epsilon} (\bar{\lambda}_{\epsilon}-x_1)}{\vert \bar{x}\vert^3}\right)<0 
		  \qquad\forall~x\in \varSigma_{\bar{\lambda}_{\epsilon}}\verb|\|B_{R_\e^1} .\] 
Notice that $\varSigma_{\bar{\lambda}_{\epsilon}}\cap\partial B_{R^{1}_{\epsilon}}=D_2\cup(D_1\cap\partial B_{R^{1}_{\epsilon}})$ with the above boundary conditions and also
		$$
		    \forall~~x\in  T_{\bar{\lambda}_{\epsilon}},\quad\Delta w^{*}_{\epsilon,\bar{\lambda}_{\epsilon}}-\frac{\alpha_{\epsilon} (\bar{\lambda}_{\epsilon}-x_1)}{\vert \bar{x}\vert^3}=0.\,
		$$
		 		  the maximum principle implies the above claim. 
		  
		 The estimate \eqref{est:w_ast} directly implies that
		  \[
		  \frac{\Delta w^{*}_{\epsilon,\bar{\lambda}_{\epsilon}}(x_1,x')-\Delta w^{*}_{\epsilon,\bar{\lambda}_{\epsilon}}(\bar{\lambda}_{\epsilon},x')}{x_1-\bar{\lambda}_{\epsilon}}\leq\frac{-\alpha_{\epsilon}}{\vert \bar{x}\vert^3} \qquad \mathrm{for~~} x=(x_1,x') \in \varSigma_{\bar{\lambda}_{\epsilon}}\setminus B_{R^{1}_{\epsilon}}.
		  \]
		  Letting $x_1\rightarrow \bar{\lambda}_{\epsilon}$, we obtain that for $\vert x\vert>R^{1}_{\epsilon}$
		  \[\frac{\partial\Delta w^{*}_{\epsilon,\bar{\lambda}_{\epsilon}} }{\partial x_1}(\bar{\lambda}_{\epsilon},x')=-2\frac{\partial\Delta v^{*}_{\epsilon} }{\partial x_1}(\bar{\lambda}_{\epsilon},x')\quad\Longrightarrow\quad \frac{\partial\Delta v^{*}_{\epsilon} }{\partial x_1}(\bar{\lambda}_{\epsilon},x')\geq\frac{\alpha_{\epsilon}}{2\vert \bar{x}\vert^3}.
		  \]
		By Taylor expansion (\ref a), we further choose $0<\delta_{\epsilon}^{1}\ll 1$ such that
		  \[\frac{\partial\Delta v^{*}_{\epsilon} }{\partial x_1}(x_1,x')\geq \frac{\alpha_{\epsilon}}{4\vert \bar{x}\vert^3}>0,\]
		  for all $\vert x\vert>R^{1}_{\epsilon}, \bar{\lambda}_{\epsilon}- \delta_{\epsilon}^{1}<x_1<\bar{\lambda}_{\epsilon}+ \delta_{\epsilon}^{1}<0$.
	\end{pf}
	\begin{lem}\label{lem3.5}
			If $\lambda_0'<\bar{\lambda}_{\epsilon}<0, w^{*}_{\epsilon,\bar{\lambda}_{\epsilon}}\not\equiv 0$, then there exist $0<\delta_{\epsilon}^{2}\ll \delta_{\epsilon}^{1}\ll 1,\,R^{2}_{\epsilon}>R^{1}_{\epsilon}$ such that
		\[
		\Delta w^{*}_{\epsilon,\lambda}>0,\quad\forall~ \vert x\vert>R^{2}_{\epsilon}, \quad x_1<\lambda,\quad\lambda \in (\bar{\lambda}_{\epsilon},\bar{\lambda}_{\epsilon}+\frac{\delta_{\epsilon}^{2}}{4}).
		\]
	\end{lem}
	\begin{pf}
	It follows from Lemma  \ref{lem3.4} that for $\lambda \in (\bar{\lambda}_{\epsilon},\bar{\lambda}_{\epsilon}+\frac{\delta_{\epsilon}^{2}}{4})$ and  $\forall~ \vert x\vert>R^{1}_{\epsilon},\quad \lambda- \frac{\delta_{\epsilon}^{1}}{2}<x_1<\lambda,\quad$

		\[
		\Delta v^{*}_{\epsilon} (x^{\lambda})>\Delta v^{*}_{\epsilon} (x)\quad \Longrightarrow	\quad \Delta w^{*}_{\epsilon,\lambda}>0.
		\]
		Here $\delta_\e^2$ is chosen such that $\lambda- \frac{\delta_{\epsilon}^{1}}{2}<\bar{\lambda}_{\epsilon}-\frac{\delta_{\epsilon}^{1}}{4}$, it suffices to consider $\vert x\vert>R^{1}_{\epsilon},  x_1<\bar{\lambda}_{\epsilon}- \frac{\delta_{\epsilon}^{1}}{4}$.
		
Notice that $\vert x\vert\sim \vert x^{\bar{\lambda}_{\epsilon}}\vert\sim\vert \bar{x}\vert\sim\vert x^{\lambda}\vert$ for all  $\vert x\vert>R^{1}_{\epsilon}$, then  
\begin{equation}\label{equ3.9}
	\left| \frac{1}{\vert  x^{\bar{\lambda}_{\epsilon}} \vert^k}-\frac{1}{\vert  x^{\lambda} \vert^k}\right|<\frac{C(\lambda-\bar{\lambda}_{\epsilon})(\vert\bar{\lambda}_{\epsilon}-x_1\vert+\vert \lambda\vert)}{\vert \bar{x}\vert^{k+2}}.
\end{equation}

By (\ref{equ3.9}) we can estimate
\[\left|\Delta w^{*}_{\epsilon,\lambda}-\Delta w^{*}_{\epsilon,\bar{\lambda}_{\epsilon}}\right|\leq
I+II+III+IV+O(\frac{1}{\vert \bar{x}\vert^4}),\]
where
\[
\begin{aligned}
I&\leq\vert b_0^{\epsilon}\vert\left|\frac{1}{\vert  x^{\bar{\lambda}_{\epsilon}} \vert}-\frac{1}{\vert  x^{\lambda} \vert}\right|,\\
II&\leq\left|\frac{b_{11}^{\epsilon}(2\lambda-x_1)^2}{\vert x^{\lambda}\vert^5} -\frac{b_{11}^{\epsilon}(2\bar{\lambda}_{\epsilon}-x_1)^2}{\vert x^{\bar{\lambda}_{\epsilon}}\vert^5}\right|,\\
III&\leq2\sum_{i=2}^{3}\left| \frac{b_{1i}^{\epsilon}(2\lambda-x_1) x_i}{\vert x^{\lambda}\vert^5}-\frac{b_{1i}^{\epsilon}(2\bar{\lambda}_{\epsilon}-x_1) x_i}{\vert x^{\bar{\lambda}_{\epsilon}}\vert^5}\right|,\\
IV&\leq \sum_{i,j=2}^{3}\mid b_{ij}^{\epsilon}x_ix_j\mid\left|\frac{1}{\vert x^{\bar{\lambda}_{\epsilon}}\vert^5}-\frac{1}{\vert  x^{\lambda} \vert^5}\right|.
\end{aligned}
\]
\end{pf}
This implies that
	\[
	\mid\Delta w^{*}_{\epsilon,\lambda}-\Delta w^{*}_{\epsilon,\bar{\lambda}_{\epsilon}}\mid<\frac{C(\lambda-\bar{\lambda}_{\epsilon})((\bar{\lambda}_{\epsilon}-x_1)+\vert \lambda\vert)}{\vert \bar{x}\vert^{3}}+O(\frac{1}{\vert \bar{x}\vert^4}).
	\]
Thus, we have
\[
\begin{aligned}
\Delta w^{*}_{\epsilon,\lambda}&\geq \Delta w^{*}_{\epsilon,\bar{\lambda}_{\epsilon}}-\mid\Delta w^{*}_{\epsilon,\lambda}-\Delta w^{*}_{\epsilon,\bar{\lambda}_{\epsilon}}\mid\\
&\geq \frac{\alpha_{\epsilon} (\bar{\lambda}_{\epsilon}-x_1)}{\vert \bar{x}\vert^3}-\frac{C(\lambda-\bar{\lambda}_{\epsilon})\left((\bar{\lambda}_{\epsilon}-x_1\right)+\vert \lambda\vert)}{\vert \bar{x}\vert^{3}}+O(\frac{1}{\vert \bar{x}\vert^4}).
\end{aligned}
\]
Since  $\bar{\lambda}_{\epsilon}-x_1> \frac{\delta_{\epsilon}^{1}}{4}$,  you can choose $\delta_{\epsilon}^{2}\ll \delta_{\epsilon}^{1},\,\,\vert \bar{x}\vert>R^{2}_{\epsilon}\gg R^{1}_{\epsilon}$,\, such that $\Delta w^{*}_{\epsilon,\lambda}>0$.
	
\begin{lem}\label{lem3.6}
	With the same assumption as in Lemma \ref{equ3.5}, there exists $0<\delta_{\epsilon}^{3}\ll\delta_{\epsilon}^{2}$ such that $\Delta w^{*}_{\epsilon,\lambda}>0,\,\,\forall~ x\in\varSigma_{\lambda},\,\, \lambda\in(\bar{\lambda}_{\epsilon},\bar{\lambda}_{\epsilon}+\delta_{\epsilon}^{3}).$
\end{lem}
\begin{pf}
	We prove it by induction. Suppose not, by Lemma \ref{lem3.5} we can choose 
	$$\lambda_{i,\epsilon}\in (\bar{\lambda}_{\epsilon},\bar{\lambda}_{\epsilon}+\frac{\delta_{\epsilon}^{2}}{4})\rightarrow \bar{\lambda}_{\epsilon},\quad x_{i,\epsilon}\in B_{R^{2}_{\epsilon}}(0)\cap\varSigma_{\lambda_{i,\epsilon}}$$
such that
	$$\Delta w^{*}_{\epsilon,\lambda_{i,\epsilon}}(x_{i,\epsilon})=
	\min\limits_{\varSigma_{\lambda_{i,\epsilon}}}\Delta w^{*}_{\epsilon,\lambda_{i,\epsilon}}\leq0.$$

 By (\ref{equc}), we know $\Delta w^{*}_{\epsilon,\bar{\lambda}_{\epsilon}}>c_{\epsilon}>0 \,\,in\,\, B_{r_{\epsilon}}(d_{\bar{\lambda}_{\epsilon}})$ with $r_{\epsilon}<\frac{\vert \bar{\lambda}_{\epsilon}\vert}{2}$. By the continuity of $\lambda$, we can choose $\delta_{\epsilon}^{3}>0$ such that 
	\[ d_{\lambda}\in B_{r_{\epsilon}/4}(d_{\bar{\lambda}_{\epsilon}}),\quad\Delta w^{*}_{\epsilon,\lambda}>\frac{c_{\epsilon}}{2}>0,\quad \forall x\in B_{r_{\epsilon}}(d_{\bar{\lambda}_{\epsilon}}).
	\]
	This imply $\vert x_{i,\epsilon}-d_{\bar{\lambda}_{\epsilon}} \vert>r_{\epsilon}$ for all $i$. 
	
	Up to a subsequence, we can assume
	 \[\lim\limits_{i\to \infty}x_{i,\epsilon}=x_{\epsilon}\not =d_{\bar{\lambda}_{\epsilon}},\quad\nabla \Delta w^{*}_{\epsilon,\bar{\lambda}_{\epsilon}}(x_{\epsilon})=0,\quad\Delta w^{*}_{\epsilon,\bar{\lambda}_{\epsilon}}(x_{\epsilon})\leq0.\]
	
	On the other hand, we can apply the proof of Lemma\,\ref{lem3.4} to show
	\[
	\Delta w^{*}_{\epsilon,\bar{\lambda}_{\epsilon}}>0,\quad x\in \varSigma_{\bar{\lambda}_{\epsilon}};\quad\Delta^2 w^{*}_{\epsilon,\bar{\lambda}_{\epsilon}}<0,\quad x\in \varSigma_{\bar{\lambda}_{\epsilon}}\verb|\|\{d_{\bar{\lambda}_{\epsilon}}\}.
	\]
	The strong maximum principle forces $x_{\epsilon}\in T_{\bar{\lambda}_{\epsilon}}$, and the Hopf lemma gives
		\[
	\frac{\partial\Delta  w^{*}_{\epsilon,\bar{\lambda}_{\epsilon}} }{\partial x_1}(x_{\epsilon})<0.
	\]
	This is a contradiction.
\end{pf}

\begin{lem}\label{f}
		With the same assumption as in Lemma \ref{equ3.5},  there exists $0<\delta_{\epsilon}^{4}\ll \delta_{\epsilon}^{3}$
		such that $ w^{*}_{\epsilon,\lambda}<0,~\forall~ x\in\varSigma_{\lambda},~
		\lambda\in(\bar{\lambda}_{\epsilon},\bar{\lambda}_{\epsilon}+\delta_{\epsilon}^{4})$.
\end{lem}
\begin{pf}
	It is known that $w^{*}_{\epsilon,\bar{\lambda}_{\epsilon}}<0$ in $\varSigma_{\bar{\lambda}_{\epsilon}}$, then $w^{*}_{\epsilon,\bar{\lambda}_{\epsilon}}<c_{\epsilon}<0 \,\,in\,\, B_{r_{\epsilon}}(d_{\bar{\lambda}_{\epsilon}}),r_{\epsilon}<\frac{\vert \bar{\lambda}_{\epsilon}\vert}{2}$. Through the continuity of $\lambda$, we can choose $\delta_{\epsilon}^{4}$ such that 
	\[ d_{\lambda}\in B_{r_{\epsilon}/2}(d_{\bar{\lambda}_{\epsilon}}),\quad w^{*}_{\epsilon,\lambda}<\frac{c_{\epsilon}}{2}<0, \quad\forall ~x\in B_{r_{\epsilon}}(d_{\bar{\lambda}_{\epsilon}}).
	\]
	By Lemma \ref{lem3.6}, we have $\Delta w^{*}_{\epsilon,\lambda}>0$ in $\varSigma_{\lambda}$, together with the boundary conditions,
	\[
	w^{*}_{\epsilon,\lambda}<0 \quad \mathrm{~on~~}\partial B_{r_{\epsilon}}(\bar{\lambda}_{\epsilon}),\qquad
	w^{*}_{\epsilon,\lambda}=0\quad\mathrm{~on~~}T_{\lambda},\qquad\liminf_{x\in\varSigma_{\lambda}, \vert x\vert\to \infty}w^{*}_{\epsilon,\lambda}(x)\leq0,	\]
 then the desired assertion follows from the strong maximum principle.
\end{pf}

	Therefore, we combine Lemmas \ref{lem3.6} and  \ref{f} to finish the proof of \textbf{Step 2}.

	\begin{step}
		$\bar{\lambda}_{\epsilon}=0$ and $v_{\epsilon}$ is radial symmetric.
	\end{step}

		Otherwise, if $\bar{\lambda}_{\epsilon}<0$, then \textbf{Step 2} implies
		 $$v^{*}_{\epsilon,\bar{\lambda}_{\epsilon}}(x)=v^{*}_{\epsilon}(x)\Longrightarrow v_{\epsilon,\bar{\lambda}_{\epsilon}}(x)=v_{\epsilon}(x).$$

This together with
		\[
		\begin{aligned}
		\Delta^2v_{\epsilon}+\epsilon\left(\frac{2}{1+\vert x\vert^2}\right)^4v_{\epsilon}&=-(v_{\epsilon})^{-7},\\
		\Delta^2v_{\epsilon,\bar{\lambda}_{\epsilon}}+
		\epsilon\left(\frac{2}{1+\vert x^{\bar{\lambda}_{\epsilon}}\vert^2}\right)^4 v_{\epsilon,\bar{\lambda}_{\epsilon}}&=-(v_{\epsilon,\bar{\lambda}_{\epsilon}})^{-7},
		\end{aligned}
		\]
	implies 
		\[\left(\frac{2}{1+\vert x\vert^2}\right)^4=
	\left(\frac{2}{1+\vert x^{\bar{\lambda}_{\epsilon}}\vert^2}\right)^4\Longrightarrow \quad \bar{\lambda}_{\epsilon}=0.\]
	This is a contradiction.

\section{Proof of Theorem 1.3}\label{sec5}
	
	\begin{lem}\label{lem4.1}
		If $u_{\epsilon}\not\equiv C_{\epsilon}$, then there exist only two antipodal critical points of $u_{\epsilon}$, where $C_{\epsilon}=(\frac{15}{16}-\epsilon)^{-\frac{1}{8}}$.
	\end{lem}
	\begin{pf}
		 We assert that there exist at least one minimal point $p_\e$ and one maximum point $q_\e$ of $u_\e$ such that  $d_{S^3}(p_\e,q_\e)=\pi$. Without loss of generality, we assume  $p_{\epsilon}=S$.
		
		If $r_{\epsilon}=d_{S^3}(p_{\epsilon},q_{\epsilon})\leq\frac{\pi}{2}$, then
		$z_{\epsilon}\in \partial B_{r_{\epsilon}}(p_{\epsilon})\cap \partial B_{r_{\epsilon}}(q_{\epsilon})\not =\emptyset$. By Theorem \ref{thm1.2}, $u_{\epsilon}$ is radial symmetric with respect to $d_{S^3}(\cdot,p_{\epsilon})$ and $d_{S^3}(\cdot,q_{\epsilon})$, then
		\[
		u_{\epsilon}(z_{\epsilon})=	u_{\epsilon}(p_{\epsilon})=	u_{\epsilon}(q_{\epsilon}).
		\]
		This contradicts the assumption that $u_{\epsilon}\not\equiv C_{\epsilon}$.
		
		If $r_{\epsilon}=d_{S^3}(p_{\epsilon},q_{\epsilon})>\frac{\pi}{2}$ and $q_{\epsilon}\not =N$, then $r'_{\epsilon}=d_{S^3}(N,q_{\epsilon})<\frac{\pi}{2}$. In $B_{r^{'}_{\epsilon}}(N)$, two cases may happen.
		
		(1) If there exists $x_{\epsilon}\in B_{r'_{\epsilon}}(N)$ such that 
		\[
		u_{\epsilon}(x_{\epsilon})=\min\limits_{B_{r'_{\epsilon}}(N)}u_{\epsilon}<u_{\epsilon}(q_{\epsilon}),
		\]
		then  $x_{\epsilon}$  is also a critical point and $d(x_{\epsilon},q_{\epsilon})<\frac{\pi}{2}$. A similar argument yields $u_{\epsilon}(q_{\epsilon})=u_{\epsilon}(x_{\epsilon})$. A contradiction!
		
		(2)\,\,$u_{\epsilon}\equiv u_{\epsilon}(q_{\epsilon})$ in $B_{r'_{\epsilon}}(N)$. By \eqref{equ1.5} we have 
		\[u_{\epsilon}(q_{\epsilon})=C_{\epsilon},\quad\,\,\mathrm{i.e.}\quad\,\, u_{\epsilon}\leq C_{\epsilon}.\]
		Integrating \eqref{equ1.5} over $S^3$ to show
		\[\int_{S^3}(C_{\epsilon}^{-8}u_{\epsilon}^8-1)u_{\epsilon}^{-7}dV_{g_{S^3}}=0.
		\]
		This implies $u_{\epsilon}\equiv C_{\epsilon}$.  A contradiction!
		
		Next we claim that there is no other critical point except for the above two critical points of $u_\e$. Otherwise, there exists a third critical point $x_{\epsilon} $ of $u_\e$, then either $d_{S^3}(x_{\epsilon},q_{\epsilon})\leq \frac{\pi}{2}$ or $d_{S^3}(x_{\epsilon},p_{\epsilon})\leq \frac{\pi}{2}$ holds. A similar argument also shows that $u_\e=C_\e$. This contradicts the assumption.

		\end{pf}
		
		\begin{lem}
		$u_{\epsilon}\equiv C_{\epsilon}.$
	\end{lem}
	\begin{pf}
		By contradiction if $u_{\epsilon}\not\equiv C_{\epsilon}$, then it follows from Lemma \ref{lem4.1} and Theorem \ref{thm1.2} that $u_{\epsilon}$ is increasing  along the great circle from the minimum point (say, $S$) to the maximum point (say, $N$). This implies that
		\[\nabla_{g_{S^3}}u_{\epsilon}\cdot\nabla_{g_{S^3}}x_{4}\geq 0.
		\]
		
		By \eqref{equ1.5} and  Lemma \ref{lem jo}, we have
		\[
		\int_{S^3}\langle\nabla(u_{\epsilon}^8\epsilon+1),\nabla x_{4}\rangle_{g_{S^3}} u_{\epsilon}^{-6}d V_{g_{S^3}}=0.
		\]
		These facts together imply that $\nabla_{g_{S^3}}u_{\epsilon}=0$, and thus $u_{\epsilon}\equiv C_{\epsilon}$ by virtue of \eqref{equ1.5}.
	\end{pf}

	{\noindent\small{\bf Acknowledgment:}  We would like to thank Professor Xuezhang Chen and Professor Xiaoping Yang for  stimulating discussions and encouragements. We  also thank Mr. Zhengyu Tao and  Mr. Nan Wu for   helpful comments on this article.

	\bibliographystyle{unsrt}

	\bigskip
	
	\noindent S. Zhang
	
	\noindent Department of Mathematics, Nanjing University, \\
	Nanjing 210093, China\\[1mm]
	Email: \textsf{mg1821015@smail.nju.edu.cn}
	
	\medskip  	
	\end{document}